\documentclass[11pt,a4paper]{article}
\usepackage{t1enc}
\usepackage[latin1]{inputenc}
\usepackage[english]{babel}
\usepackage{harvard}
\newcommand{\N}{\mbox{\rm I\hspace{-0.31ex}N}}
\newcommand{\Z}{\mbox{\rm Z\hspace{-0.8ex}Z}}

\newcommand{\C}{\mbox{\rm\hspace*{0.6ex}\rule[0.2ex]{0.08ex}{1.2ex}\hspace{-0.7ex}C}}
\newcommand{\R}{\mbox{\rm I\hspace{-0.33ex}R}}

\newcommand{\beq}{\begin{equation}}
\newcommand{\eeq}{\end{equation}}
\newcommand{\beqarr}{\begin{eqnarray}}
\newcommand{\eeqarr}{\end{eqnarray}}
\newcommand{\beqa}{\begin{eqnarray*}}
\newcommand{\eeqa}{\end{eqnarray*}}
\newtheorem{theorem}{Theorem}

\begin{document}
\thispagestyle{empty}
\begin{center}
\section*{Introducing Groups into Quantum Theory \\(1926 -- 1930)}
{\em Erhard Scholz, Wuppertal}  \\[0.5pt]
\end{center}

\begin{abstract}
In the second half of the  1920s, physicists and mathematicians introduced group theoretic methods into the recently invented ``new''  quantum mechanics. Group representations  turned out to be a highly useful tool in spectroscopy and in giving quantum mechanical explanations of  chemical bonds. H. Weyl  explored the possibilities of a group theoretic approach towards quantization.  In his second version of  a gauge theory for electromagnetism, he even started to build a bridge between quantum theoretic symmetries  and  differential geometry.    Until the early 1930s, an active group of  young  quantum physicists and mathematicians contributed to this new challenging field. But around the turn to the 1930s, opposition against the new methods  in physics  grew. This article focusses on the work of those  physicists and mathematicians who introduced group theoretic methods into quantum physics.
\end{abstract}

\subsection*{Introduction}
In the middle of the 1920s,    understanding of the representations of Lie groups
and understanding  of  the quantum mechanical  structure
of matter   made great advances, almost
simultaneously. 
Certain members of both disciplines saw  the
potential for building   new and deep connections
between mathematics and theoretical physics.
Thus a cooperative development highly consequential  for
theoretical physics began in the second half of the 1920s, with the main protagonists being   W. Heisenberg, E. Wigner, F. London,  W. Heitler and, to a lesser degree, P.A.M. Dirac on the  one side, H. Weyl, J. von Neumann,  and B.L. van der Waerden on the other. The first 
introduction and  use of the new method in theoretical
physics met soon with opposition (``group pest''). But  it
turned out to be successful in the long run, and to be just the first wave
of a process of restructuring  mathematical concepts and techniques
in the theory of the basic structures of matter. After an intermediate period of  about two decades with a slow and nearly unnoticed continuation of work in this direction,   another wave of using group-theoretical methods in physics gained momentum  in the second half of the century. This development has
recently attracted interest from the side of  history and
philosophy of science.\footnote{\cite{Mehra/Rechenberg:VI,Gavroglu/Simoes,%
Karachalios:Diss,Brading:Symmetries} }
It should be quite as interesting from the point of view of  the history of mathematics, because it established  broad and consequential semantical relations for an important field of modern mathematics.

The following article explores  the first wave of introduction
of  new mathematical methods into  quantum physics and
chemistry. It starts
with the early realization of the usefulness of group theoretic
methods for the study of spectroscopy and chemical bonds, and
stops short of the consolidation of what was  achieved in the
first wave in three textbooks on the
subject published in the early 1930s, \cite[$^2$1931]{Weyl:GQM}, \cite{Wigner:Gruppen} 
 and \cite{Waerden:GQM}, which  have now become classics of the field. Unlike the
other two, Weyl's book  had an earlier  first edition at the end of the 1920. It therefore enters
 the period of investigation of our investigation.

  This article is   a first step into this interdisplinary terrain
 from the side of history of mathematics. It relies heavily on the solid background laid out
by T. Hawkins'  study \cite{Hawkins:LieGroups} and H. Rechenberg's chapter on group theory and quantum mechanics in
 \cite[VI 1,  chaps. III.4, III.5]{Mehra/Rechenberg:VI}.

\section*{1. Heisenberg and Wigner}

Shortly after the invention of the new quantum mechanics,   P.A.M. Dirac, W. Heisenberg,  and E. Wigner  started to consider   consequences of symmetry  in
multi-particle systems for the structure of energy terms in
atomic spectra.\footnote{For the emergence of matrix, wave, and ``q-number''  mechanics  see, among others,
\cite{Hendry:BohrPauli,Beller:Qdialogue,Pais:Inward,Rechenberg:QM,%
Cassidy:Heisenberg,Kragh:Dirac,Moore:Schroedinger}. A multi-volume encyclopedic report is 
 \cite[vols, II, III,  IV, V]{Mehra/Rechenberg}. A  six-page compression of the crucial period 1923 -- 1926 can be found in the introduction to volume VI  
\cite[VI 1, xxv--xxxi]{Mehra/Rechenberg:VI}. For a splendid bibliography  see
 \cite[VI 2, 1253--1439]{Mehra/Rechenberg:VI};  indexes of the whole series  at the end of the same volume VI 2. } 
  Dirac studied the role of antisymmetry  in multi-electron systems in summer 1926. Important as that was for the growing understanding of  quantum mechanics, it did not employ group theory beyond the distinction of the signum of permutations. Group theoretic questions  proper started to be addressed by Heisenberg  and Wigner in late 1926 and early 1927.

The newly established paradigm of quantum mechanics demanded to characterize  a (quantum) physical system, at the time typically an electron system in the shell of an atom or of a molecule, by a set of Hermitian (or more generally, symmetric) operators, one for any observable quantity of the system, in a state space $\cal S$ assumed to be a Hilberts space in order to have sufficient symbolical  structure. In   Schr\"odinger's perspective, $\cal S$  was  viewed as a space of complex ``wave'' functions. Then the tool of differential operators could be used.\footnote{Questions  how the function space was to be completed, or how domains of the operators should be understood and, perhaps, extended, were generously  neglected by the early quantum physicists. Such questions were first  addressed  by J. von Neumann in the later 1920s and at the turn to the 1930s.}
 Most important was the operator characterizing the energy of the system (or a constitutive part of it, like an electron in a multiparticle system), the {\em Hamilton operator}
$H$. Other operators could characterize linear momenta $P_i$ or coordinatized  spatial  positions $Q_i$ ($1 \leq i \leq 3$), rotational  (orbital) momenta $L_i$, the square of the total momentum $L^2$, and, a little later, the {\em spin} $J$  of a particle (considered to express the ``particle's'' proper rotation) etc.. 

For an  
 atom, the eigenspaces of the Hamilton operator $H$ could  characterize the stationary states of a system of electrons, or of an outward electron, depending on the situation. The eigenvalues $E_1, E_2 , \ldots$  of  $H$  represented the energy values obtained in these states. Often such eigenstates turned out to be degenerate,  i.e., they belonged to an eigenvalue of multiplicity $> 1$. This was the case for  atoms or molecules  with rotational  symmetry.  Of course, spectroscopy did not allow to measure the energy of each eigenstate  directly. Only differences between two energy values, say $E_1$ and $E_2$, were  observable by the frequency $\nu$ of the radiation emitted during the transition of an electron from one  energy state to the other, 
\[    h \nu = E_1 - E_2 \; . \] 

In  early 1925,  Pauli conjectured that bound electron states in a molecule have  an intrinsic two-valuedness and that electrons obey  an {\em exclusion principle}  forbidding different electrons (a littler later also other ``fermions'') to occupy the same state of a system. Later in the year, S. Goudsmit and E. Uhlenbeck established the hypothesis of   electron {\em spin} which they assumed to arise from a   ``proper rotation'' of the electron. Different empirical evidence indicated that this intrinsic spin was  quantized with respect to any specified spatial direction in exactly two possible states $u$ and $d$ (spin ``up'' and spin ``down''). Early in 1927, W. Pauli mathematized the idea by a  spin state space $\C^2$ extending the complex phase of the
 Schr\"odinger wave function $\psi (x)$ \cite{Pauli:Spin}. In group theoretic language, which was not yet in Pauli's mind in early 1927, he implicitly worked inside the natural representation of $SU_2$, the covering group of the spatial rotations $SO_3$. He  proposed to describe a spinning particle by a two valued wave function $\tilde{\psi } = ( \psi_1 , \psi_2 ) $, later called a {\em Pauli spinor}.\footnote{Pauli drew upon the symbolic ressources of the 
Klein-Sommerfeld theory of the spinning top, which contained the natural representation of $SU_3$ implicitly. For a review of the understanding of the rise of spin see \cite{Waerden:Spin} or \cite[chap. VI.4]{Mehra/Rechenberg:IV}.}
 It could be constructed from  Schr\"odinger  wave functions by forming  (tensor) products with the complex two-dimensional space characterizing the complex superpositions of the two possible pure spin states  $\C ^2 \cong <u, d> $  (here $< \quad > $ denotes the linear span).  The total wave function of a collection of $n$ electrons was  expressed formally as a ``product'' (in later terminology  as an element of the $n$-fold tensor product  $\bigotimes^n {\cal S}$).
 In summer 1926 P.A.M. Dirac  realized  that Pauli's exclusion principle implied  that  
multi-electron (more generally fermion) states had to  be represented by {\em alternating products} \cite{Dirac:Antisymm}.\footnote{Cf. \cite{Kragh:Dirac}} 

\subsubsection*{An ad-hoc usage of permutations (W. Heisenberg) }
Already before Pauli's mathematization of spin was known, Heisenberg started to consider the consequences of the new phenomenon for multi-electron systems. In June 1926 he  submitted his first  paper on this topic to  {\em Zeitschrift f\"ur
Physik} \cite{Heisenberg:Permutationen}. He looked for 
reasons for the separation of energy terms in the spectrum of
higher atoms into different subsets between which  apparently no
exchange of electrons took place ({\em term systems without
intercombination}). Such an effect could be seen by ``missing''  lines when one compared the observed spectral lines with  the combinatorics of all the arising energy levels in a higher atom. Heisenberg guessed that  the interaction of the
orbital magnetic momentum of electrons (i.e., the magnetic
momentum resulting from what was left from Bohr's electron orbits
in the new quantum mechanics) with the  still hypothetical  spin
might play a crucial role for this phenomenon \cite{Heisenberg:Permutationen}. 

In a
second part of the paper, submitted in December 1926, he
continued to  explore  the hypothesis further. He proposed the view that the
distinction of term systems might result from  a kind of
``resonance phenomenon'' between the spin states of the different electrons and, perhaps, their orbital momenta.
He made clear  that here the word
``resonance'' was    not to be understood in the sense of
classical physics, but  as an expression of  a  physical intuition of the ``more subtle interplay of the electrons in an atom''
 \cite[556, 578]{Heisenberg:PermutationenII}. Thus Heisenberg's ``quantum mechanical resonances'' referred to spin coupling effects for which at that time no  adequate mathematical representation was known.\footnote{In early quantum chemistry the term  ``resonance'' was  used in a comparable metaphorical way; see \cite{Mosini:Resonance}.}
He therefore  looked for new tools to deal with them and hoped to find them in  the theory of permutation groups.

In his investigation, 
Heisenberg studied   states  of  $n$-electron systems in an atom or molecule. Abstracting at first from spin, he started from 
 $n$      eigenfunctions  $l, m, \ldots p$ 
(Heisenberg's notation)  of the  Hamilton operator, which described possible states of single electrons without spin,  possible degeneracies  included.
As usual he  described  a composite system   by a kind of noncommutative product of the eigenfunctions. He considered the result as a state of the ``unperturbed'' composite system, while the spin coupling (``resonance'') had to be taken into account as a perturbation due to the ``more subtle interplay of the electrons''.
Because  electrons are indistinguishable, he  concluded:
\begin{quote}
In the unperturbed case, the eigenfunction of the total system
can be written as product of all functions of the single
electrons, e.g., $l_1 m_2 \ldots p_n$. The unperturbed problem is
$n!$-fold degenerate, because a permutation of the electrons
leads to equal energy values of the total system.
\cite[557]{Heisenberg:PermutationenII}
\end{quote}

For an element $u$ of the (tensor) product space, written by our author as
 $u = l_1 m_2 \ldots p_n$ with an  index $1 \leq i \leq n$ for the
different electrons, Heisenberg considered the result of an
electron permutation $ S \in {\cal S}_n$, the symmetric group of $n$ elements,
 and wrote it as
\[   S u =  l_{S (1)} m_{S (2)} \ldots p_{S (n)} . \]

 If we  denote the state space of a single electron by $V = <l, m, \ldots , p>$,  $ dim\, V = n$, the ($n !$-fold degenerate) total  state space of the quotation above
corresponds to the span of  vectors arising from  permutation
of the components of any one product state $u$.\footnote{We may prefer to distinguish Heisenberg's  basic state vectors by  a lower index $i$,  $\psi_1 =l, \psi_2 = m, \ldots, \psi_n = p$, and to characterize the bijection between states and  electrons by adding an upper index $j$, $\psi_i^{(j)}$ ($1 \leq, i, j \leq n$). Then it is advisable to 
order the tensor product according to  electron indexes, 
$\psi_{i_1}^{(1)} \otimes \psi_{i_2}^{(2)} \otimes \ldots \psi_{i_n}^{(n)}  $ (comparable to  Wigner's notation, see below). That  makes the upper (electron) index redundant, and the lower (state) index $i$  encodes the different possibilities for bijections completely. Because Heisenberg ordered according to states and used the electron indexes to indicate the bijection between electrons and individual states, his permutation $S$  operated on the state vectors  of the (``our'') tensor product  $V^{(n)}$ by inversion   $S^{-1} =: \sigma $, i.e. from the right:   $(\psi_1 \otimes \psi_2 \otimes \ldots \otimes \psi_n ). \sigma  = \psi_{\sigma (1)} \otimes \ldots \otimes \psi_{\sigma (n)} = \psi_{S^{-1} (1)} \otimes \ldots \otimes 
\psi_{S^{-1} (n)}$.  As this detail has no consequences for the orthogonality questions, we  follow Heisenberg's description  in the sequel without further retranslations. } 
We  we want to denote
it here as $V^{(n)}$ 
\[  V^{(n)}  :=  \langle S u \, | \, S \in {\cal S}_n  \rangle \subset \otimes{}^n V  . \]
 $V^{(n)}$ was constructed to  characterize  the state space of an ``unperturbed'' system of $n$ electrons distributed according to Pauli's principle (i.e., mapped bijectively) on  the $n$ states $m, l, \ldots, p$. Without spin the energy was totally degenerate (all eigenvalues identical), while the consideration of spin split it up into different ``non-combining'' terms. The physical model of the electron system had to account for the impossibility of  transitions of electrons between the respective states or subspaces. Mathematically the question was whether the corresponding vectors (wave functions) or subspaces in Hilbert space were   orthogonal.
 
Heisenberg  looked for a decomposition of $V^{(n)}$ into ``noncombining'' (orthogonal) subsystems if spin resonance was considered as a kind of perturbation. As we will see in a moment, he had good arguments that orthogonality of subspaces should not be affected by the spin perturbation. Its basic structure could thus be analyzed already on the level of the unperturbed system without spin.

In order to address this question, Heisenberg considered a cyclic subgroup of ${\cal S}_n$ generated by a ``substitution'' (permutation) $S$ of highest possible order $\nu$, and  an orbit  in $V^{(n)}$ of an eigenstate $u$ under such a subgroup. He then formed different superpositions   of the elements of such an orbit.  For a permutation $S$ of order $\nu$ he choose coefficients formed by  powers of a primitive $\nu $-th root of unity $\omega, \; \omega ^{\nu} = 1$, in the following way:
\beqa
U_0 &=& \frac{1}{\sqrt{ \nu } }(u + Su + S^2 u +\ldots S^{\nu -1}u) \\
U_1 &=& \frac{1}{\sqrt{ \nu } }(u + \omega Su + \omega ^2 S^2 u +\ldots \omega ^{\nu -1}S^{\nu -1}u) \\
& \ldots & \\
U_{\nu -1} &=& \frac{1}{\sqrt{ \nu } }(u + \omega^{\nu -1} Su + \omega ^{2(\nu -1)} S^2 u +\ldots \omega ^{(\nu -1)^2}S^{\nu -1}u) .
\eeqa

These linear combinations were  formed in   analogy to the construction of the roots of resolvents in  the  theory of algebraic equations. In fact, Heisenberg referred to a textbook of higher algebra, a fifty year old German translation of a classical  book by Serret \cite{Serret:German}, which had been written originally in 1866 (third edition), as one of the first books containing a passage on the recently revived theory of E. Galois.\footnote{\cite[110ff.]{Kiernan:Galois}} 
For dimensional reasons ($\nu < n!$) there were elements $w = T u, \; T \in {\cal S}_n$,
 of the defining basis of  $ V^{(n)}$ (Heisenberg: ``eigenfunctions'') which were linearly independent of the $U_0,  \ldots , 
U_{\nu -1}$. They  lead to analogously formed linear superpositions $W_0, \ldots, W_{\nu -1}$. He applied  the same procedure, step by step,   until the whole space  $V^{(n)}$ was  spanned by elements of such a  form: $U_0, \ldots, U_{\nu -1 }, W_0, \ldots, W_{\nu -1 }\, \ldots  $.\footnote{Cf. \cite[489ff.]{Mehra/Rechenberg:VI}.}

Now, Heisenberg
 collected all functions $U_j, W_j , \ldots$ starting with the  same exponent $j$ of the unitary root $\omega $ into one collection,
\[ \Gamma _{\omega ^{j}} := \{ U_j , W_j, \ldots  \} \,  , \]
and proposed that the corresponding subspaces could be taken as  symbolical representatives for the different  term systems. He  argued that  the span of  $\Gamma _{\omega ^{j}}$ and  $\Gamma _{\omega ^{k}}$ ought to be orthogonal (for different $j$ and $k$) 
\beq \int \bar{f_j} g_k = 0  \; , \qquad 
 \qquad f_j \in \Gamma _{\omega ^{j}} , \; g_k \in \Gamma _{\omega ^{k}} , \qquad  j \neq k \; .\eeq
His argument for this claim depended  crucially on  an invariance argument of the transition integral under any permutation:\begin{quote}
If under the integral (\ldots)  the electron numbers are somehow permuted, the value of the integral cannot change. \cite[559]{Heisenberg:PermutationenII}\end{quote}
The physical context of the calculation demanded such an invariance. Although  Heisenberg's construction of the ``term systems'' $\Gamma _{\omega ^{j}}$ did not ensure  such an invariance, his argument held for   similar constructions in  which  the invariance condition was satisfied.\footnote{We will  see in a moment (equation 
 (\ref{standard-rep-S3})) that  Heisenbergs $\Gamma _{\omega ^{j}}$, respectively their linear spans, are no invariant subspaces under the full permutation group  . Heisenberg's own argument shows that therefore his model was   physically unreliable.  Wigner's approach solved the problem. It was different to Heisenberg's, contrary to what the latter believed. }
 The form of his argument was close to  one used in early Galois theory (``as the whole constellation does not depend on the choice of the ordering of the roots of the equation, \ldots  such and such inference can be drawn \ldots'') and may have been prompted by the latter.

 Heisenberg agreed with Dirac that an ``eigenfunction'' of the total system should be antisymmetric under permutation of the electrons. It seemed  impossible, at the moment, to draw consequences of this
 postulate.\footnote{A structural answer to this question was given  later by Weyl and a  more pragmatic one  by  von Neumann and Wigner, see below.}
 On the other hand, he plausibly assumed  that any {\em perturbation} of transition probabilities, arising from spin coupling, should be symmetric under transposition of two electrons. That was sufficient, in his context, to  show that the decomposition of the total space of $n$ electrons $V^{(n)}$ into orthogonal subspaces was not affected by  spin resonance. Thus, so he concluded, the subspaces spanned by the $\Gamma _{\omega ^{j}}$  ought to characterize  the decompositions of energy terms into non-combining partial systems  {\em including} spin \cite[559]{Heisenberg:PermutationenII}.\footnote{I thank an anonymous referee for  having made me aware of this  important passage in Heisenberg's argument.}  
 Although the argument  did not work in his  own ad-hoc construction, it would become  important (and correct) once it was transferred to a decomposition into truely invariant subspaces. 

All in all, Heisenberg's paper gave an  inventive treatment of the term system problem, although it must have apppeared  surprising for   mathematical readers of the time (like  J. von Neumann or H. Weyl). For the 
construction of  non-combining term systems, Heisenberg relied on a rather  
old-fashioned algebraic background  \cite{Serret:German}.  Neither H. Weber's textbook \cite{Weber:Algebra} nor any other more recent algebraic text was  even
mentioned. Such a neglection of more recent methods may not   necessarily  be  of great disadvantage for a new application of mathematics by a physicist. But in this case, the neglection of younger algebraic developments  included the  methods of representation theory of finite groups, which dealt with structures much closer to Heisenberg's problem than   algebraic equation theory. In his  first step into the new terrain, Heisenberg  had to rely on formal expressions originally introduced in a  completely different context.  Thus his hypothesis for the identification of 
non-combining term systems by his $\Gamma $-collections was quite daring and would surely  have led to difficulties, had it been used in future investigations without major modifications.

From hindsight it is easy to see that Heisenberg's decomposition did   not lead to   irreducible representations of the permutation group. Worse than that, Heisenberg's 
hypothetical ``non-combining term systems''   $\Gamma _{\omega ^{j}}$ {\em were not even invariant subspaces}  under the full permutation group. His construction made sure  that a  subspace $\Gamma _{\omega ^{j}}$  is an eigenspace with eigenvalue $\omega ^j$ of the cyclic subgroup generated by the permutation $S$. But this does not hold  for other permutations. Already for $n=3$,
 $\omega = e^{\frac{2 \pi i }{3}}$ and any   
 3-cycle $S$, e.g. $ S = (123)$, a transposition $T$ with $TST=S^2$, e.g. $T=(12)$, maps  $U_1 \in \Gamma _{\omega }$ to $U_2 \in \Gamma _{\omega ^{2}}$, 
\beq  \label{standard-rep-S3} SU_1 = \omega U_1 \, ,   \quad STU_1 = T S^2 U_1 = \omega ^2 T U_1 \; . \eeq
In fact,  the linear spans of $\{ U_1, U_2  \}$ and $\{ W_1, W_2  \}$, in Heisenberg's notation, are  copies of the two-dimensional irreducible representation of ${\cal S}_3$.\footnote{Cf. \cite[8ff.]{Fulton/Harris}.}
In other words, the irreducible spaces are  transversal to the subspaces offered by Heisenberg as ``non-combining term systems''. 
But before such discrepancies could start to irritate other contributors to the  program, Heisenberg's method was outdated by an approach to the problem proposed by his colleague  E. Wigner.  

So it was good news, and even better ones than Heisenberg knew, that  he could  refer  to  Wigner's investigations already  in a footnote added in proof to his December paper.  He  erroneously believed that his approach agreed with Wigner's    \cite[561, footnote (1)]{Heisenberg:PermutationenII}. In fact, a rash view could support this belief, as in the case of 3-electrons, e.g.  a Lithium atom, both methods led to   equal numbers and dimensions  of the respective term systems: two 
one-dimensional  term systems  (symmetric and  antisymmetric) and two equivalent 
two-dimensional term systems (standard representation in Wigner's approach), $ 6 = 1 + 1+ 2+ 2$.    But while Wigner characterized the non-combining term systems by subspaces which actually were  irreducible subrepresentations, we have seen that Heisenberg's decomposition  was different, even in this  case.  

In the end, it appears as a lucky sequence of events that 
  Wigner's  papers threw   new light on the question so fast.  His approach superseded Heisenberg's group theoretically ad-hoc method, before  the latter could lead into a dead end. Wigner's papers  opened  the  path towards an introduction of  group representation into the study of  multi-particle systems and established a sound mathematical frame into which Heisenberg's perturbation calculation could be integrated without contradictions.\footnote{In the literature on history of quantum mechanics this essential difference between Heisenberg's and Wigner's approaches is often passed over in silence; cf. e.g.,    \cite[489ff.]{Mehra/Rechenberg:VI}.}
 
\subsubsection*{Turn towards group representations (E. Wigner)}
Eugene Wigner had studied chemical engineering at Budapest and
Berlin (TH) during the years 1920 to 1925 and had gained access to
the physical community organized around the colloquia of the {\em
Deutsche Physikalische Gesellschaft} and the local {\em
Kaiser-Wilhelm Institutes}.\footnote{For the following
passage on Wigner compare \cite{Chayut:Wigner} and
\cite{Mackey:Wigner}.}
 After he had finished his diploma
degree, he went back to Budapest and worked as a chemical
engineer in a leather tannery  (his father's craft), but he
continued to read the {\em Zeitschrift f\"ur Physik} with the
interest of an {\em aficionado}. Thus he was well informed about
the breakthroughs in quantum mechanics, achieved during 1925. He
immediately accepted the chance to go back to Berlin, when he  was
invited by  Karl Weissenberg to become his assistant   at the
{\em Kaiser-Wilhelm Institute} for fibre research. Weissenberg
himself had studied applied mathematics with R. von Mises and had
then  turned towards condensed matter physics. He needed support
in his X-ray  investigations of crystal structures. At
Weissenberg's suggestion, Wigner started to read the group
theoretic parts of Weber's textbook \cite{Weber:Algebra} and
 to explore the symmetry characters of crystals in the new
setting.\footnote{See \cite{Chayut:Wigner} and Wigner's
autobiographical report in \cite[105]{Wigner:Auto}.}
 Because of this  interest in actual X-ray crystallography,
he  was much better acquainted with group theory than Heisenberg in
1926. 

 In late 1926, Wigner  started to study the question of how
$n$-particle systems can be built from $n$ given, pairwise
different, single particle states $\psi _1, \ldots, \psi_j, \ldots
\psi _n$, initially without considering spin effects.   Like Heisenberg,  he wanted to know how the $n$-particle
state space decomposes  under permutations of the electrons.
Each  electron was  (in the stationary case) identified mathematically  by its hypothetical
``space coordinates''  $r_i  =(x_i, y_i,
z_i) \in \R^3$, where $i$ served as an  index to characterize
different electrons.

In his first paper on the topic \cite{Wigner:Terme_I}, submitted
on November 12, 1926, he  considered a  product of $n$
``eigenfunctions'' $\psi _1, \ldots, \psi _n$. Any state
$\psi _k$ can be ``occupied'' by any  (the
 $i$-th) electron, which was denoted by Wigner by $ \psi _k (r_i)$. He then
considered permutation states of the form
\[ \psi_{\sigma 1} (r_1) \psi_{\sigma 2} (r_2) \ldots \psi_{\sigma n} (r_n) =: v_\sigma , \]
where $\sigma $ is any permutation of $n$ elements, (the notation
$v_\sigma $ is ours). Thus Wigner studied essentially the same
subspace  $V^{(n)}$ of the $n$-fold tensor product of $V = <\psi _1, \ldots,
\psi _n>$ as Heisenberg. 
In his first paper he considered only the special
case $n=3$  and calculated  the decomposition of  $V^{(3)}$ into
irreducible components under permutations ``by hand''. No wonder, that he found Dirac's symmetric and antisymmetric representations
among them and  in addition two 2-dimensional ``systems''.\footnote{The regular representation of ${\cal S}_3$  (cf.  next footnote), $R_3 \cong V^{(3)}$,  decomposes into the trivial representation $U$, the antisymmetric representation $U'$ 
(both 1-dimensional) and two copies of the twodimensional irreducible subspace $S_2 := \{ (z_1, z_2,z_3) \, | \, z_1 + z_2 + z_3 = 0\}$ of the natural representation on $\C^3$ arising from permutations of the basis vectors: $R_3 = U \oplus U' \oplus S_2 \oplus S_2$.   }
 He concluded similar to Heisenberg:
\begin{quote}
The additional systems are all degenerate, this degeneration is
such that it cannot be broken by any perturbation symmetric in the
single particles which are assumed to be equivalent.
\cite[34]{Wigner:Terme_I}
\end{quote}

The state space  $V^{(3)}$  was  spanned by vectors $v_{\sigma }$ identified  by  permutations $\sigma \in {\cal S}_3$. The operation of  ${\cal S}_3$ on $V^{(3)}$ was multiplication of permutations (in Wigner's case from the left), just like in the {\em regular representation}.\footnote{The {\em regular
representation} $R_G$ of a finite group $G$ is given by the
operation of $G$ on the group algebra $\C[G]:= \{  \sum_{h}  z_h
h | z_h \in \C \}$ (summation of $h$ over $G$) by operation from
the left. It contains all finite dimensional  irreducible
representations of $G$. More precisely, in each representation of the
symmetric group of $n$ elements each irreducible component $X$
appears in the regular representation with multiplicity $dim \,
X$. Cf. \cite{Fulton/Harris} or any other book on representation theory.}
  In this way Wigner hit, at first unknowingly,  upon the problem
of a decomposition of the regular representation of the symmetric
group ${\cal S}_3$. His approach to the problem made it apparent that, more generally,  $V^{(n)}$ was by its very  construction just another version of the  regular representation of the symmetric group.  It had been studied  by Frobenius,
Schur, Burnside, Young  and others in their   works   on the
representation theory of finite groups.\footnote{See \cite{Hawkins:Characters_I,Hawkins:Characters_II} and the overview in \cite[373--384]{Hawkins:LieGroups}.}  

When Wigner discussed this question  with J. von Neumann, a good
friend of his since their common school days at Budapest, his friend immediately recognized what Wigner was doing from a
mathematical point of view and explained  the  problem in terms
of a decomposition of the regular  representation. Thus Wigner started the second part of his contribution (submitted
November 26, 1926) with a general observation which introduced
the representation theory of the symmetric group. Noting the
rising calculational complexity, when one wanted to extend the
results from $n=3$ to higher cases, he remarked:
\begin{quote}
There is a well  prepared mathematical theory, however, which one
can use here, the theory of  transformation groups isomorphic to
the symmetric group (\ldots ), which has been founded at the end
of the last century by Frobenius and has been elaborated later by
W. Burnside and J. (sic!) Schur, among others. J. von Neumann was
so kind to make me aware of these works, and predicted the
general result correctly, after I told him the result for the
case $n = 3$. \cite[43]{Wigner:Terme_II}
\end{quote}
Therefore Wigner considered it worthwhile  introducing the basic
facts of the representation theory of the symmetric group to the
readers of the {\em Zeitschrift f\"ur Physik}.\footnote{For a more recent introduction to the subject,   see  \cite{Sternberg:Groups}.}
  In particular, he explained in  his
article how on can calculate the dimension $N_{(\lambda )}$ of   a
representation of ${\cal S}_n$ characterized by a partition
$(\lambda ): = (\lambda _1, \ldots, \lambda _k)$ of
$n$,\footnote{The dimension of  $N_{(\lambda )}$  is the quotient of $ n !$ by the product of all ``hook lenghts'' of the corresponding Young diagram. For  details  see  \cite[89ff.]{Sternberg:Groups}.}

\[  n = \lambda _1 + \lambda _2 + \ldots + \lambda _k, \;\;\; \lambda _i \geq \lambda _{i+1}  .\]

After Wigner became aware of the decomposition of the regular
representation, he  could adapt Heisenberg's perturbation argument for spin coupling to the modified context:
\begin{quote}
In a system with $n$ equal mass points, between which initially
there is no exchange of energy, each eigenvalue is $n!$ degenerate
(if the corresponding state does not contain equivalent orbits).
If one creates  an exchange of energy, each eigenvalue splits into
several. \cite[44]{Wigner:Terme_II}
\end{quote}
He proposed to calculate the degeneracy of the corresponding term by the dimension $N_{(\lambda )}$ as above. The basic structure for the  splitting  of energy terms in an atom with $n$ (peripheral)  electrons, which had been translated by Heisenberg into  the problem of decomposing $V^{(n)}$ into minimal  invariant  subspaces, was  now elucidated  by applying  standard methods of representation theory  for the symmetric group. To Wigner and von Neumann 
this turn may have appeared   like some kind of 
 ``pre-established harmony'' between physics and mathematics, stipulated  in the contemporary G\"ottingen milieu of mathematics and mathematical physics. For other participants it may have looked more like a kind of magic of mathematical symbolism.

 On the other hand, many questions were still open. Among them  most importantly the question  which of   the irreducible representations of the permutation group on the space of Schr\"odinger wave functions were compatible with the Pauli-Dirac principle  of antisymmetry for the total (Pauli-) wave function.  In order to address this question, the spin phenomenon and its relation to rotational symmetries had to be understood better.

\section*{2. Wigner and von Neumann}
Early in 1927, Wigner made considerable advances. He  enriched
the study of invariance by including rotations
of the state space of electrons in an  outer atomic shell.
 In his third paper in spectroscopy, he started to derive the basic structural data of
  spectroscopic terms  from the rotational symmetry of the electron
state spaces \cite{Wigner:Terme_III}.\footnote{Recieved  May 5,
1927.} Already in the introduction to the paper he stated:
\begin{quote}
The simple form of the Schr\"odinger differential equation allows
us  to apply  certain group methods, more precisely,
representation theory. These methods have the advantage that  by
their help one gets results nearly without calculation, which do
not only hold {\em exactly}  for the one-particle problem
(hydrogen atom), but also for arbitrarily complex systems. The
disadvantage of the method is that it does not allow us to derive
approximative formulas.   In this way it is possible to explain a
large part of our qualitative spectroscopical experience.
\cite[53]{Wigner:Terme_III}
\end{quote}

\subsubsection*{Representations of the rotation group}
Again it was J. von Neumann who advised Wigner what to read in order
to understand the representation theory of the special orthogonal
group $SO_3 $, in particular the  recent papers by I. Schur and
H. Weyl \cite{Schur:1924,Weyl:1924[61]}.\footnote{See \cite[63,
fn. (1)]{Wigner:Terme_III}.} Thus Wigner discussed, among others,
the irreducible representations of the rotations in the plane,
$SO_2$, which are  (complex) 1-dimensional. They are characterized by
an integer parameter $m$, such that any plane rotation $\delta
_{\alpha }$ by an angle $\alpha $ has the representation as the
(one by one)  ``matrix''  $e^{i m \alpha }$.  Let us denote, for
 brevity, this representation of the plane rotation
group as $d^m$. Then, of course, the representation matrix of the
rotation  $\delta _{\alpha }$ is the $1\times 1$ matrix
\[   d^m ( \delta _{\alpha })  = e^{i m \alpha }; \]
in other words, the representation of the rotation by the angle $\alpha $ has the eigenvalue $e^{i m \alpha }$.

Wigner then introduced the  $(2l+1)$-dimensional representations
of  $SO_3$ (of highest weight $l \in \N_0$), which we denote here
as ${\cal D}^{l}$, according to present conventions, and
indicated how to calculate the representation matrices
\[ D^l (A) = ( D^l_{j k}(\alpha , \beta ,\gamma ))_{1\leq j,k\leq 2l+1} \]
for any rotation $A \in SO_3$, characterized by its three Euler
angles $\alpha , \beta ,\gamma $ \cite[68ff.]{Wigner:Terme_III}.
Moreover, he discussed the decomposition of ${\cal D}^{l}$ under
restriction to the subgroup $SO_2$ of rotations about the
$z$-axis into $2l+1$ one-dimensional subspaces. This leads to 
representations $d^{m}$ in our notation above, where $m$ may assume
the $2l+1$ pairwise different values
\[  -l \leq m \leq l . \]

That fitted structurally  so well with the observed
classification of spectra and their discrete parameters, the  quantum
numbers,  that Wigner could immediately proceed to a
spectroscopical interpretation of these representation theoretic
quantities. The highest weight $l$ could be identified with the {\em azimuthal quantum number} of  the 
Bohr-Sommerfeld theory   \cite[71]{Wigner:Terme_III} (later often called {\em orbital angular momentum} quantum  number).\footnote{In spectroscopy, an alphabetical code is used for   $l$: $S$ for $l=0$, $P$ for $l=1$, $D $ for $l=2$ etc.. }
 Moreover, the weight $m$ of the specified abelian subgroup  $SO_2$  appeared  as a group theoretic characterization of the {\em magnetic quantum number} of the electron. The latter had been introduced in order to explain the split of spectral lines (indexed by the {\em principal} quantum number $n$ of the so-called Balmer-series and by $l$) into different terms (``multipletts'') under the influence of a  strong magnetic field, the so-called normal {\em Zeeman effect}.\footnote{With a magnetic field in direction of the observation, P. Zeeman had observed such an effect in 1896, while perpendicular to the field  a ``third'' (undisplaced) line appeared. H.A. Lorentz had explained it a year later in terms of a classic theory of the electron in the magnetic field, cf. \cite[161]{Rechenberg:QM},   \cite{Darrigol:Electrodynamics} or   \cite[76f., 268ff.]{Pais:Inward}.}
A similar effect had been observed under the influence of a homogeneous electric field
  ({\em Stark effect}).\footnote{The Stark effect had been observed in 1913.}
Thus the basic features of the dynamics of  the electron were apparently closely related to the basic parameters of   representations of the symmetry group of its orbit.

After a short discussion  of the fact that transitions of electrons occurred
in nature only between neighbouring  azimuthal (orbital angular momentum) quantum numbers $l$, corresponding to a change  $\triangle l = \pm
1$, Wigner turned to the consequences of the introduction of a
homogeneous electric field:
\begin{quote}
By means of an electric field along the $Z$-axis the substitution
group of our differential equation is diminished (verkleinert).
Thus we have to proceed  [as above] and reduce the
three-dimensional rotation group  to a collection of
representations of the two-dimensional group (about the $Z$-axis).
\cite[72]{Wigner:Terme_III}
\end{quote}
As a result, under the influence of an external homogeneous
field, a term with azimuthal quantum number $l$  splits into
$2l+1$ lines, indexed by the magnetic quantum number $m$.\footnote{In this context (Stark effect), Wigner called $m$ the ``electric quantum number'' \cite[73]{Wigner:Terme_III}.} 

 For atoms with more than one electron involved in radiation processes, the situation was, of course, much more complicated. Here Wigner could only vaguely indicate, how the representation of the rotation group and of permutations might work together to form the the total state space of an $n$-electron system and how they determine the combined quantum numbers \cite[77f.]{Wigner:Terme_III}.

\subsubsection*{The spin group $SU_3$}
For a  detailed investigation,  a more
subtle study of the interplay between rotational symmetry, its
relation to spin  properties, and the exchange symmetries
(permutations) of multi-particle systems became necessary. At
almost the same time as Wigner's paper on  rotational symmetries,
Pauli submitted his path-breaking proposal to mathematize
Uhlenbeck's and Goud\-smit's hypothesis of an intrinsic ``spin'' of
the electron  by the use of ``two-component'' wave functions
\cite{Pauli:Spin}.\footnote{Received May 8, 1927, by {\em
Physikalische Zeitrschrift}, three days after the submission of
Wigner's paper \cite{Wigner:Terme_III}.} 
Charles G. Darwin stepped in with a series of papers on the ``electron as a 
vector wave''.\footnote{\cite{Darwin:electron_I,Darwin:electron_II}}
That made it possible for
Wigner  to extend  the investigations   of  symmetries  to  spin effects.

 For such studies  von Neumann's advice became even more important than before.
 The publications discussed above were  written by E. Wigner when
he was still an assistant for theoretical chemistry at the
technical university Berlin. In spring 1927 he moved  to
G\"ottingen  for one year, as an assistant of Hilbert's. At that
time, Hilbert  suffered strongly from pernicious anemia and was
nearly  inaccessible to his new assistant. Nevertheless, Wigner
came into close contact with other young physicists working at
G\"ottingen, among them in particular L. Northeim, P. Jordan, and
W. Heitler.  Moreover, von Neumann visited  G\"ottingen regularly
\cite{Mehra:Wigner}. Thus  there were good conditions for  Wigner
and von Neumann to establish the basic representation theoretic
features of atomic spectra, including  spin effects, during late
1927 and the first half of 1928, simultaneously with H. Weyl's
work on the same topic and independently of it.

Between December 1927 and June 1928, E. Wigner and von Neumann
submitted a series of three  papers on spectra and the ``quantum
mechanics of the spinning electron (Drehelektron)''  to the {\em
Zeitschrift f\"ur Physik}.\footnote{Dates of reception: December
28, 1927; March 2, 1928; June 19, 1928.}
 As Wigner  later reported,  he wrote the papers after intense  discussions
 with his colleague and friend whom he therefore considered to be
a coauthor \cite[496]{Mehra/Rechenberg:VI}. In this series, the
authors emphasized the conceptual  role of representation theory
for quantum mechanics in an explicit and programmatic manner and parallelized it to the invariance method of general relativity.
\begin{quote}
\ldots It may not be idle to call  the strong heuristical value
(Sp\"urkraft)  to attention, which dwells in these and similar
principles of symmetry, i.e. invariance, in the search for the
laws  of nature: In our case it  will lead us, in a unique
and compelling way,  from Pauli's qualitative picture of the
spinning electron to the regularities of the atomic spectra. That
is similar to the  general theory of relativity, where an
invariance principle made it possible to unveil the universal
laws of nature. \cite[92]{Neumann/Wigner:I}
\end{quote}

In their paper, Wigner and von Neumann took up Pauli's characterization of spin by   a (commutative) product of a
Schr\"odinger wave function
\[  \psi (x), \;\;\; x = (x_1, x_2, x_3)  \in \R^3, \]
and a complex function  $\zeta (s)$ depending on  variable in a discrete two-point ``internal'' spin space, $s \in
\{\pm 1  \} $. The combined function
 \beq \varphi (x,s) = \psi
(x) \zeta (s)   \eeq
 had been introduced by as  Pauli  as (spin-) wave function. The
 dependence on $s$ could just as well be written in  index form
\[ \varphi_s (x) :=  \varphi (x,s) , \;\;\; \mbox{with} \;\; s \in \{ \pm 1 \}   . \]
Then the Pauli wave function was given by two components,
\[ \tilde{\varphi}  (x) := (\varphi_{-1} (x),  \varphi_{1} (x) ) \, ,   \]
and $\tilde{\varphi} $ could be considered as a modified wave function (on $\R^3$) with
values in $\C^2$, a ``hyperfunction'' in Wigner's terminology (later called a {\em Pauli spinor field} on $\R^3$).

For  an $n$-particle system the wave function acquired the form
\beq  \label{Pauli-wf}   \tilde{\varphi}  (x_1, \ldots , x_n) := (\varphi_{s_1 \ldots s_n}(x_1, \ldots, x_n) ) \;, \;\;\; x_j \in \R^3, \; s_j \in \{ \pm 1 \}\; . \eeq
Then the values of $\tilde{\varphi} $ were in $\C^{2n}$ \cite[94]{Neumann/Wigner:I}

Wigner and von Neumann studied how to express the operation of
the rotation group $SO_3$ on the Pauli wave-functions by a unitary
operator. They introduced an explicit expression for the
complexified version $ \tilde{A}$ of a rotation $A  = A(\alpha
,\beta ,\gamma ) $ given in terms of the Euler angles $\alpha ,
\beta , \gamma $
 \cite[98]{Neumann/Wigner:I},
\beq \label{complex-rotation}
\tilde{A} := \left(  \begin{array}{rr}  e^{- i \frac{\alpha }{2}}& 0 \\   0 &  e^{i \frac{\alpha }{2}} \end{array}  \right)
\left(  \begin{array}{rr}
 \cos \frac{\beta }{2}&  \sin \frac{\beta }{2}   \\  - \sin \frac{\beta }{2}& \cos \frac{\beta }{2}  \end{array}  \right)
\left(  \begin{array}{rr}
  e^{-i \frac{\gamma }{2}}& 0 \\
0 &  e^{i \frac{\gamma }{2}}  \end{array}  \right)
 .\eeq

\[  A \mapsto \tilde{A} \in SU_2,  \]
such that a rotation $A^{-1} \in SO_3$   operated on the wave-functions by
\beq \label{opSO3} \varphi (x) \mapsto \tilde{A} \varphi (A^{-1} x) . \eeq

That agreed  well with what Pauli had done; but while Pauli  had
made use of the complex description of the spinning top, well
known in the Sommerfeld school, Wigner and von Neumann  embedded
the  formula into a representation theoretic perspective. In
particular they referred to the second paper of Weyl's great
series on the representation theory of the classical Lie groups
\cite[98, footnote]{Neumann/Wigner:I}. Here Weyl had discussed
the universal  coverings of the   special orthogonal groups
(later to be called {\em spin groups}), had proved the  full
reducibility and derived the characters and dimensions of all
irreducible representations \cite{Weyl:Darstellung}.\footnote{See
\cite{Hawkins:LieGroups}.}
 Von Neumann and Wigner stated clearly  that they  needed
only certain  aspects of the general theory.\footnote{``Of
course, much less than Weyl's deep rooted results are necessary
for our present goals.'' \cite[98, footnote]{Neumann/Wigner:I} }
But they made quite clear that now one had to take into account
``two-valued'' representations of the $SO_3$, in addition to the
(one-valued) ones studied by Wigner in his last paper (called
above  ${\cal D}^l$, $l \in \N_0$).
 That gave an additional series which  will be  denoted here by ${\cal D}^{\frac{k}{2}}$ ($dim({\cal
D}^{\frac{k}{2}}) = k+1$), $k$ odd,
according to more recent conventions.\footnote{Cf. \cite[181ff.]{Sternberg:Groups}.}

For the goal of their paper, they considered the most basic
two-valued representation, in fact a local inverse of the covering
map
\[ SU_2 \longrightarrow SO_3 ,\]
 given by equation (\ref{complex-rotation}) up to sign.
Then ${\cal D}^{\frac{1}{2}}$ was given by  the standard
representation of the covering group $SU_2$; more precisely
\[   {\cal  D}^{\frac{1}{2}}A = \pm \tilde{A} . \]
  In the perspective of  their paper, this representation arose naturally from
the operation of  $SO_3$ on the 1-particle state as described
in equation (\ref{opSO3}). It  was  essential to find the consequences
for the $n$-particle state.

They indicated how to find the  matrix expressions of a
representation matrix ${\cal D}^{\frac{k}{2}}A$ for a rotation $A \in
SO_3$, characterized by its Euler angles $\alpha ,\beta, \gamma
$, in analogy to  Wigner's formulas in the classical (one-valued)
case. In doing so,  they relied on Weyl's result and stated that
for each dimension $n  \in \N$ there exists exactly one
representation of $SO_3$ (or its universal cover) indexed by $j:=
\frac{n-1}{2}$. In the sequel we use the slightly more recent unifying notation for the two series: 
\beq \label{representation SU2} {\cal D}^j =   {\cal D}^{(j,0)} , \;\;\; \mbox{of dimension}\; n
= 2j + 1  , \;\;\;  j \in \{0, \frac{1}{2}, 1, \frac{3}{2}, 2,
\ldots  \}
 \eeq
Here $n$ odd (respectively $j $ integer valued) corresponds to
one-valued representations, and $n$  even ($j$ half-integer) to
``two-valued'' representations of the orthogonal group.

With the machinery of representation theory at their disposal, it
was clear how to proceed  to the description of  the $n$-particle
states described by $n$-fold tensor products. They ended the
first paper of the series with an observation on how to decompose
the tensor product spaces into irreducible components:
\begin{quote}
In the applications it will be important to know the irreducible
representations of the rotation group in $\{ ^n a^{({\cal
R})}_{s,t} \}$  [Wigner/von Neumann's symbol for $\otimes^n {\cal
D}^{\frac{1}{2}}$, E.S. ]; that is easily achieved, as its trace
is additively composed from the traces of the former.
\cite[108]{Neumann/Wigner:I}
\end{quote}
They gave an explicit result, described verbally, but without any  ambiguity.
 Written  in more recent symbolism, it was
\beq  \otimes^ n {\cal D}^{\frac{1}{2}} = {\cal D}^{\frac{n}{2}} \oplus (n-1) {\cal D}^{\frac{n-2}{2}} \oplus \frac{n}{2} (n-3) {\cal D}^{\frac{n-4}{2}} \oplus \ldots \; . \eeq

\subsubsection*{Permutations,  spin, and anomalous Zeeman effect}
In the second paper of their series, Wigner and von Neumann
combined the rotational and spin symmetries with the permutation
aspect from which Wigner had started. Wigner's basic physical
intuition was that in atomic spectroscopy the energy operator $H$
will be composed,
\[ H = H_1 + H_2, \]
 by  a part $H_1$ resulting from the spatial motion of the electron only
 (the motion of the ``center of gravity'' of the electron, as he
said) and the ensuing gross effect of the electromagnetic
interaction  with the field of the atomic core. The second
part, $H_2$, should model other aspects, most important among
them the electron spin \cite[133]{Neumann/Wigner:II}. Thus one
could start from the eigenvalue problem of the ``spin-less'' wave
function $\psi $ (Schr\"odinger wave function),
\[   H_1 \psi = \lambda \psi \; , \]
and  refine  the result by passing to
the ``hyperfunctions'' $\varphi$ including spin (i.e., the Pauli spinors).

For the investigation of  symmetry  properties with respect to
permutations, it was therefore natural to  distinguish  different
types of operations for a permutation $\alpha \in {\cal S}_n$, an
operation $P$  on {\em space variables only} 
and an operation $O$ on  {\em both}, spin and space variables 
 ($P_{\alpha }$ and $O_{\alpha }$ in Wigner's notation):
 \beqa P_{\alpha }^{-1} \varphi (x_1, \ldots,
x_n;s_1, \ldots, s_n) & := &
 \varphi(x_{\alpha _1}, \ldots, x_{\alpha _n};s_1, \ldots, s_n ) \\
O_{\alpha }^{-1} \varphi (x_1, \ldots, x_n;s_1, \ldots, s_n) & := &
 \varphi(x_{\alpha _1}, \ldots, x_{\alpha _n};s_{\alpha _1}, \ldots, s_{\alpha _n} ) .
\eeqa
The operation $Q$ of permutations on {\em spin variables only}  could be constructed from these \cite[133]{Neumann/Wigner:II} by
\[ Q_{\alpha } := P_{\alpha }^{-1} O_{\alpha } . \]
Obviously ``spin-less'' wave functions transformed under $P_{\alpha
}$, while the transformation $O_{\alpha }$ of ``hyperfunctions''
could be built from $P$ and $Q$,  $ O_{\alpha }  = P_{\alpha } Q_{\alpha }$.

Wigner then considered a slow continuous change  from an energy state in which the
spin contribution could be neglected ($H = H_1$) to one, in
which this was no longer the case
\cite[133]{Neumann/Wigner:II}.
He  made the following observation: 

While the  original state with $H=H_1$ is invariant under $O$ and $P$, an increasing  spin  perturbation $H_2$ may reduce  the original  symmetry to $O$ only. In this case,  the
formerly  irreducible subspaces for $H_1$ are decomposed into
smaller irreducible components of $H_1 + H_2$.

That was  a convincing  group theoretic view of the split of spectral terms by a  perturbation bringing spin differences into the game. Empirically such a phenomenon had been observed long ago in the {\em anomalous Zeeman effect}: If  a weak magnetic field was switched on, spectral lines belonging to the same magnetic number $m$ could split into different terms.\footnote{The ``anomalous Zeeman effect'' had been observed by A.A. Michelson and T. Preston in 1898, and could not be explained in the Bohr-Sommerfeld theory of the atom; cf \cite[161f.]{Rechenberg:QM} or \cite{Pais:Inward}.}

But it was still to clarify how to deal with the antisymmetry principle for the total wave function of an $n$-electron system. 
According to Dirac ``\ldots only those states
occur in nature, the eigenfunctions of which are antisymmetric''
 \cite[133]{Neumann/Wigner:II}.  
 Wigner and von Neumann therefore continued with the study of the irreducible representations
 of the symmetric group ${\cal S}_n$ in the antisymmetric part of
the total ``hyperfunction'' representation, i.e., in
\[ \wedge^n \tilde{V}  \subset \otimes^n \tilde{V },  \]
where $\tilde{V}$  denotes a state space of  single-particle  ``hyper-functions''
(Pauli-spinor fields). Of course, such irreducible antisymmetric representations are one-dimensional, and the question was, under which conditions such antisymmetric representations in the ``hyperfunction'' space could be derived from  an irreducible representation of the spin-free wave functions. 
To simplify language,  we denote the representation of ${\cal S}_n$ in $V^{(n)}$ corresponding to 
 a partition $(\lambda )=(\lambda
_1, \ldots, \lambda _k)$   by $V^{(n)}_{(\lambda )}$.

If one
starts from a degenerate energy term with multiplicity $m$ of the spin-less Schr\"odinger equation of an
 $n$-electron system
\beq \label{degnerate}   H_1 \psi = E_0 \psi \; ,  \eeq
 one can form a basis of $m \, 2^n$ corresponding ``hyperfunctions'',  by allowing for the combinatorics of possible spin values for the $n$ constituents.
If  analogously $m$ denotes the dimension of  an irreducible representation  $V^{(n)}_{(\lambda )}$ like above,  the $m \, 2^n$-dimensional space of spin extended hyperfunctions may be called   $\tilde{V}^{(n)}_{(\lambda )}$. Obviously it    forms   an invariant subspace of $  \otimes^n \tilde{V }$ (under permutations). Our authors now looked for  irreducible components of  $\tilde{V}^{(n)}_{(\lambda )}$, and in particular one-dimensional antisymmetric ones.

Using a result of A. Speiser's book on group theory \cite{Speiser:Gruppen}, they came to the conclusion 
 that a partition   $(\lambda )$ allows to form a (non-trivial, one-dimensional) antisymmetric extension in $\tilde{V}^{(n)}_{(\lambda )}$, if   and only if  $(\lambda )$  is of the form 
\beq  \label{2-2-1-1}  (\lambda ) = (2, 2, \ldots , 2, 1, 1, \ldots , 1) \, . \eeq

That was an important result for the group theoretical 
  program in spectroscopy. It showed clearly, {\em  why (and under which conditions) irreducible representations of the symmetric group could  characterize a term system} of higher atoms.

Still the question had to be answered, in how many fine structure terms a   spectral line of an $n$-electron system, corresponding to an azimuthal (orbital momentum) quantum number $l$ and partition $(\lambda )$,  could split. Thus 
Wigner and von Neumann finally  studied 
 the  combinatorical possibilities,  by which  the total magnetic quantum number
 $m = m_1 + \ldots m_n$  of such a system could be  built  from the
 quantum numbers $m_j$ of the individual electrons  and which  effects  could be expected from switching on a spin perturbations  $H_2$.    They came to the conclusion that the   momentum (including spin)    of an $n$-electron  system in such a state  can be characterized by a  (integer or half-integer)  value  $j$, called {\em internal quantum number}, with
\[  |\frac{n-2z}{2} -l  |   \leq j \leq \frac{n-2z}{2} +l \; \]
(with  difference 1 betweeen two values of $j$). For each $j$ the total magnetic  momentum including spin  $\tilde{m}$ then may acquire values in $-j \leq \tilde{m} \leq j$.
 The number
 $t$ of different values for $\tilde{m}$, i.e., the number of possible terms into which the $n$-electron state $(\lambda)$ with azimuthal quantum number $l$   could split, was then, according to \citeasnoun[140--143]{Neumann/Wigner:II}: 
\[  t = \mbox{min}\;  \left\{  \begin{array}{l} n- 2z + 1 \\
2l + 1 \qquad \qquad \mbox{}
\end{array}  \right.     \]

This result agreed beautifully with empirical  findings and with the rules derived in other theoretical approaches.\footnote{Like Hund's ``Aufbauprinzip'' \cite[140]{Neumann/Wigner:II}.}
 Wigner was proud about what he had achieved  cooperatively  with von Neumann:
\begin{quote}
Thus the, probably, most important qualitative spectroscopical
rule has been derived. Independent of the immense effectiveness
(Leistungsf\"ahigkeit) of quantum mechanics (\ldots), one will be
surprised that all this was ``plucked out of the air'', as one
might say (da\ss{} alles, wie man sagt ``durch die Luft'' ging ), i.e., without taking into account the special form of
the Hamiltonian function, only on the basis of symmetry
assumptions  and of Pauli's qualitative idea.
\cite[143]{Neumann/Wigner:II}
\end{quote}
Although definite values of the  energy differences could not be derived by  group theoretic methods alone, Wigner's and von Neumann's approach  gave a convincing  explanation for the  splitting of 
a spectral line under a  magnetic field  (Zeeman effect) of any kind into ``multiplett''  terms of the   fine structure. 

\section*{3. London and Heitler}
  In quantum chemistry, representations of
permutation groups made their first appearance  about the same time as they did in spectroscopy. The topic was opened up by a joint
publication of  two young physicists, Walter Heitler and Fritz
London,  who had come to Z\"urich  on Rockefeller grants in 1926
(F. London), respectively  1927 (W. Heitler), to work with E.
Schr\"odinger.\footnote{\cite{Gavroglu:London}}
 While a closer scientific cooperation with their professor turned out to be
 more difficult than expected, they used the opportunity to
exchange and develop ideas with each other. In June 1927 they
submitted a paper on the quantum mechanical explanation of so-called
covalent bonds (those due to valence electron pairs), which arose from an idea of W. Heitler. It  soon
was considered as the entry point for quantum mechanical model
building in chemistry \cite{Heitler/London:1927}. According to L.
Pauling,  one of the great figures  of the first generation in
quantum chemistry,  Heitler's and London's paper can be
considered as
\begin{quote}
\ldots the greatest single contribution to the clarification of
the chemist's conception which has been made since G. Lewis's
suggestion in 1916 that the chemical bond between two atoms
consists of a pair of electrons held jointly by two atoms
\cite[340]{Pauling/Wilson:1935} (quoted from
\cite[542]{Mehra/Rechenberg:VI}).
\end{quote}

 The  story of this invention leads deep into the history of quantum theory and of
 chemistry  and is covered as such in the respective historical
literature.\footnote{See
\cite{Gavroglu/Simoes,Karachalios:2000,Karachalios:Diss,Nye:Chemistry,%
Simoes:Cambridge}
and \cite[540ff.]{Mehra/Rechenberg:VI}.} We want to concentrate
here on a   specific  aspect, 
which is at the center of our investigation of the use of modern
mathematical methods in physical chemistry: the contexts, reasons
and mode for the appearance and use of group theoretic methods.
Such methods were first applied   in two papers by W. Heitler, published
in 1928 \cite{Heitler:Gruppen_I,Heitler:Gruppen_II}. They built upon  a joint paper with F. London, published
during their common summer in  Z\"urich \cite{Heitler/London:1927}.

In their joint paper, Heitler and London started from an investigation
of two hydrogen atoms and their electrons, initially modelled
separately, at a distance  $d = \infty $ between the nuclei,
by identical Schr\"odinger functions with energy eigenvalue
$E_0$.    Using a perturbative  approach, they studied what
happened to the electrons and their added energies when  the
atomic distance $d$ was reduced. They showed the existence of  two
solutions, $\psi_1$ and $\psi_2$  for the combined system, with
respective total energies $E_1$ and $E_2$, and interpreted the
energy difference
\[ \triangle E _i := E_i - 2 E_0 \; , \;\; i = 1, 2 ,    \]
as a kind of {\em exchange energy} of the electrons.\footnote{The
quantum physical idea behind this terminology was the following:
If one   joined  two probability ``clouds'' about two nuclei  to
one (of the combined system) some kind of ``exchange'' of
particles between  two ``partial clouds'' related to the nuclei,
although   fused to represent    one state, seemed now possible
(i.e., had positive probability).  The language of ``exchange
energy'' has to be taken, again, as a classical metaphor for a
quantum effect. For a more detailed discussion see
\cite[380f.]{Schweber:Slater}. } With their choice of sign,
negative exchange energy expressed that the compound system had a
lower energy state than the two single systems. Moreover, the
exchange energies were dependent on the distance parameter $d$.
Their analysis showed that,  with $d$ increasing from a little
above 0 to some value $d_1$, $E_1$ fell to a minimum, rising
again for increasing $d$  from $d_1$ to $\infty $, while $E_2$
fell monotonously for $d > 0$ with increasing $d$ ($ d
\rightarrow \infty$). Thus $\psi_1$  represented a bound state
for $ d = d_1$,  while $\psi_2$ characterized  a repulsive force
for any value of the atomic distance (the van der Waals repulsion
between the two hydrogen atoms)\cite[460]{Heitler/London:1927}.

A continuation of the calculation for two helium atoms, each
containing two electrons,  showed that  only the case of a repulsive
interaction could be obtained, if  electron spin and the
Pauli exclusion principle were taken into account (i.e., if both
electrons of one atom were assumed to be in different spin
states). In this sense, the ``exchange energy'' of Heitler and
London appeared as an effect of spin coupling and was positive in this case. 
It explained why helium did not form two-atomic
molecules and behaved as noble gas. The principles of
non-relativistic quantum mechanics seemed to open the possibility
of understanding the {\em structure} (graph-like combinatorics of
atomic ``valences'') and the {\em quantity} (energies)  of {\em
chemical bonds}.

\subsubsection*{Heitler's theory of valence  bonds}
In summer 1928, E. Schr\"odinger went  from Z\"urich to Berlin,
as a successor on M. Planck's chair; in  October  F. London
joined him there as an assistant.  W. Heitler, whose Rockefeller
grant had run out more or less at the same time, accepted an
offer from Max Born to become an assistant at G\"ottingen. There
he got to know E. Wigner whose group theoretic works he had 
started to read with great interest  when in
Z\"urich.\footnote{\cite[VI.1, 502, 547]{Mehra/Rechenberg:VI}} Now
Heitler explored what  the representation theory of the symmetric
group could achieve for the determination of quantum mechanical
bond states.

 Already in   January 28, 1928, he  submitted his first article on the topic
\cite{Heitler:Gruppen_I}. His goal was to extend the approach of
his joint work with London to    ``higher'' molecules. For the time being, 
 that did not mean more than   two-atomic molecules with $n >
2$ outer electrons. He stated his methodological preferences
clearly at the beginning of the paper:
\begin{quote}
Among all methods, the group theoretic is the one which
definitely achieves most for the multi-particle problem: it was
brought in by E. Wigner [Heitler referred to
\cite{Wigner:Terme_II,Wigner:Terme_III}, E.S.] to achieve a
qualitative overview of  all existing terms.
\cite[836]{Heitler:Gruppen_I}
\end{quote}

Heitler came to the conclusion that {\em already at large
distances} the exchange forces between valence electrons of {\em
opposite spin} resulted in a reduction and even a relative
minimum of bond energy, which expressed an attractive force
between the two atoms. Here he defined {\em valence electrons} as
such electrons of quantum numbers $(l, m)$ in the outer
``shell''\footnote{``Outer shell'' now referred  to electrons of
highest azimuthal (orbital momentum) quantum number $l$ with respect to its spherical symmetry
${\cal D}^l$  in the atom, and with a compatible magnetic quantum number
$m$ ($-2l \leq m \leq 2l$).}
 which had no partner  of equal quantum numbers $l$,  $m$ with 
 opposite spin in the same atom.  Heitler hinted at certain
restrictions of his approach:
\begin{quote}
We still have to warn of an overestimation of the implications
(Tragweite) of our results in two respects. The simple formulas
for the interaction energy \ldots can only be considered as a very
rough approximation, because the perturbative calculation neglects
several points and holds  only for large distances. Secondly, the
``exchange molecules'' considered by us represent only a part of
the chemical molecules. although of the most prominent and most
stable ones ($N_2, \, O_2, \, NH_3, \, CH_4$ etc.). A large part
of the homopolar compounds, however,  relies on perturbations of a
different kind \ldots .\footnote{Heitler referred to the
neglection of ``polarization'' which he estimated  for $H_2$  to
be about 25 \% and guessed that it should be much higher for
higher molecules.} \cite[837]{Heitler:Gruppen_II}
\end{quote}

Thus Heitler distinguished clearly between different kinds of chemical bonds only some of which could be explained, in his opinion, by spin coupling accessible to  group theoretic methods. He called them {\em exchange molecules}.
We have  to keep this  in mind when we look at the extension
 of Heitler's theory of valence bonds from
a more structural,  mathematical point of view (e.g., by   Weyl)
and its reception by physicists and chemists.

Here, Heitler investigated  two electron systems $A$ and $B$, each of
which consisted of $n$ (valence) electrons, initially without interaction.
All in all, he studied  a system of $2n$ electrons. Following
Wigner, he  characterized  a term system by an irreducible
representation of the permutation group of $2n$ elements ${\cal
S}_{2n}$. Let us call it $R$.

Under the assumption of no interaction, $R$ could  also be
considered as a  representation of each of the $n$ electrons $A$
and $B$ and thus of two subgroups isomorphic to ${\cal S}_n$, let
us say $R_A$ and $R_B$. The latter were no longer
irreducible. Thus Heitler studied the decomposition of $R$ into
subspaces which were  simultanously irreducible    in $R_A$ and
in $R_B$. This work was  facilitated by the assumption (unproved
but considered as self-evident by Heitler) that the 
Pauli  principle implies that 
\begin{quote}
\ldots the representations appearing in nature [are] those which contain
only 2 and 1 in their partition \cite[846]{Heitler:Gruppen_I}.\footnote{This condition was proved  a little later by
Wigner in his  joint work with von Neumann, as we have seen. It  may have been orally communicated knowledge in  G\"ottingen  already in winter 1927/28.}

\end{quote}

He concluded that only  those  representations could appear, in 
which for both partial systems $A$ and $B$ the respective $n$
valence electrons are characterized by  a completely ``antisymmetric term
system''  and have antiparallel spin
\cite[848]{Heitler:Gruppen_I}. On this basis he was able to give an
approximative calculation of the exchange energies.

This result established a
 quantum mechanical explanation of certain   non-ionic bonds  which could not be explained in
terms of Coulomb forces. Traditionally, chemists had used
{\em valence dashes} to represent such molecules. In 1916, G. Lewis had
proposed a qualitative interpretation of a valence dash as a pair
of electrons shared by two atoms. But the underlying
physical forces  remained a  mystery. Now it seemed promising to
look for an explanation of such ``valences'' by   the pairing of
electrons with opposite spin, but otherwise equal quantum
numbers. Heitler's proposal was thus to investigate the range of
the hypothesis that  {\em spin coupling of valence electron
pairs} lay at the base of  molecule
formation.

In  a second article on the topic, submitted September 13, 1928,
Heitler extended his investigations to molecules with more than 2
atoms \cite{Heitler:Gruppen_II}.
Here Heitler was less cautious than  in January. He now  described the result of his first article as having
established   a  ``complete equivalence'' of the quantum
mechanical explanation of  homopolar chemical bonds for
two-atomic molecules and the traditional explanation  of chemical
valences by electron pairs (Lewis). He introduced an
integral expression $J_{Q}$ derived by Heisenberg for the
exchange energy  between two systems $Q$, constituted by the
partial systems  $A$ and $B$ \cite{Heisenberg:Ferromagnetismus}, and resumed:
\begin{quote}
Each such exchange energy $J_Q$ between two atoms can be
interpreted  as a {\em valence bond}
symbolically denoted by a  valence dash (Valenzstrich).  Nearly
all typical and stable two-atomic molecules of chemistry rely on
such an exchange bond; and vice versa: if the valence theory
permits the existence of a two-atomic molecule then it is
possible  quantum mechanically. \cite[805, emphasis in
original]{Heitler:Gruppen_II}
\end{quote}

Although his theory did not predict new or different effects in
comparison to classical chemical knowledge,  it claimed to explain the
empirical knowledge of valence bonds  structurally, for the case
of {\em two-atomic} molecules. Moreover, it should lead to  a
quantitative determination of bond energies, even if  only in
the sense of  a rough, first estimation (see quotation above).

\subsubsection*{Other approaches}
Competing approaches to the quantum mechanics of chemical bonds were developed 
by F. Hund and a little later by L. Pauling, R. Mulliken, and
others. They  shed doubt on the range of
Heitler's and London's theory and on its   quantitative reliability. They did not rely on the exchange energy
of spin coupling, but  concentrated on the spatial distribution of the Schr\"odinger function. During the next decade it
turned out that for more complicated molecules Heitler's method
led to unrealistic predictions. The      alternative approaches
were necessary, even on the   structural level, to achieve a
satisfactory  agreement with experimental knowledge.

In summer 1928 these consequences were not yet clear, although
chemists like Mulliken and Pauling already  thought along
different lines.\footnote{See
\cite{Gavroglu/Simoes,Nye:Chemistry} or
\cite[552ff.]{Mehra/Rechenberg:VI}.} For a short while Max
Delbr\"uck who    became well known for his later researches  on
the molecular basis of genetics   considered Heitler's and London's
approach  worth following. He studied perturbative formulas for
the determination of exchange energies based on group theoretical
methods \cite{Delbrueck:Gruppen}. Thus Heitler could see his
position strengthened and contributed to further explorations of
his method in \cite{Heitler:Gruppen_II}. Here he posed the 
fundamental   question as to the {\em existence} of  multi-atomic
molecules, on the basis of exchange energies of valence pairs of
electrons.

 This type of question was highly interesting from a mathematical point of
 view, but may have appeared 
useless  to most  chemists.   Heitler considered his investigation as
nothing more than a ``preliminary study (Vor\-stud\-ie)''. In the
course of it, he came to admit that in the calculations of
exchange energies,  it might happen that permutations of more
than two electrons contribute essentially to the interaction. 
That  had already been conjectured by F. London. Heitler  remarked
that, in his  opinion, bonds which rely on such higher exchanges
could not be considered as ``valence bonds in the sense of
Lewis''. They would constitute a different type of bond. Nevertheless he
thought it justified  to study,  how far one could come with
valence bonds  proper (``in the sense of Lewis'')
\cite[815]{Heitler:Gruppen_II}. At the time, he  still hoped that
chain molecules of organic chemistry and lattice structures might
belong to ``our bond category'' \cite[806]{Heitler:Gruppen_II}.

This hope did not come true. During the 1930s, L. Pauling's and
R. Mulliken's approach of constructing  ``molecular orbitals'',
i.e., Schr\"odinger functions of  multi-electron systems about a
complex of atoms (molecular core), built much less on 
structural principles such as permutations. They drew upon previously
unformalized chemical knowledge on hypothetical spatial
constellations of the atoms for the modelling of  Schr\"odinger
functions of a system of electrons. The striking successes of
this approach  turned out to be crucial for the acceptance of
quantum mechanics among chemists \cite{Gavroglu/Simoes}. It
became the core mathematical technique during the next few
decades for a  fruitful elaboration of quantum mechanical models
for more complicated molecules, in particular in organic
chemistry.\footnote{Up to our days, it continues to be  the basis  for the semi-classical approximations used as building blocks for  the computer simulations of molecular structures, cf. \cite{Le Bris/Lions}.}

 \section*{3. Weyl at the backstage}
Taking the  results of Wigner,  von Neumann, and
Heitler  into account, it might look as if not much was  left for
Hermann Weyl when he entered the field. But such an impression would be completely wrong; Weyl 
took up a whole range of  questions  pertaining to the challenging new field and  entered into 
 second phase of  active involvement in mathematical physics between
 1927 and 1931. This second phase  was a  natural follow up to
his first phase of activity in theoretical physics between
1917 and 1923, in which he had made crucial contributions to
general relativity,   unified field theory, and
cosmology.\footnote{ See
\cite{Skuli:Diss,Coleman/Korte:DMV,Scholz:DMV,Mackey:Kiel,Speiser:Kiel}.}
 When he entered the terrain of quantum mechanics, he was  particularly
 interested in the role of  group representation  and
contributed to  the introduction of gauge methods into the
quantum physical setting.

The background of  Weyl's intervention in the field was one of
the surprising conjunctions in the history of science, which
turned out to be  tremendously fruitful. During the years 1925/26
the M\"unchen-G\"ottingen-Cop\-en\-hag\-en group of Heisenberg,  Born,
 Jordan, and Pauli,  closely communicating with  Bohr, invented
quantum mechanics;  Schr\"odinger,  at that time working at
Z\"urich, complemented it  with  his ``wave mechanics'',  
P.A.M. Dirac, in Cambridge,  developed his perspective of ``q-numbers'' (a formal operator symbolism, particularly well adapted to the physicists way of thinking) and crowned the whole development by an
overarching view -- called  ``transformation theory'' by
physicists. 

At the beginning of this period, in April 1925, Weyl
had just finished his great work on  the representation theory of
classical (Lie-) groups.\footnote{Weyl delivered the three parts
of the series \cite{Weyl:Darstellung} in January, February, and
April 1925. For this part of the story see
\cite{Hawkins:LieGroups,Borel:Essays,Slodowy:Weyl}. } For him, it
was not only the attraction of  the fascinatingly rich
mathematical structures of  covering groups,  decomposition of
representations into  irreducible spaces, calculation of
characters, classification of root systems, weight  vectors, and reflection groups
etc., which made  him turn towards this work, but rather its
intriguing interplay   with  conceptual questions lying at
the basis of physical theory building. Weyl had met classical
groups and Cartan's classification of their infinitesimal
versions (Lie algebras) on two occasions during his first phase of
active involvement in mathematical physics. He  found them to be
crucial for  answering  two  questions in this context:
\begin{itemize}
\item[--] Why are tensors such a good and, in fact, universal tool in
general relativity and, more general, in differential geometry?
\item[--] What are group theoretic reasons for the ``pythagorean''
(Weyl's terminology for what later was called semi-Riemannian) nature of the metric in  general relativity?
\end{itemize}

The {\em first question} was answered by Weyl in 1925 with the
insight, and its proof, that all irreducible representations of
the general linear group $GL_n\R$ can be constructed  as invariant 
subspaces of  tensor powers of the underlying standard
representation (for differential geometry,  $V = T_pM \cong \R^n$,  the tangent space  at a point $p$
to the underlying manifold $M$). In this sense, tensors and tensor spaces were universal objects for the representation of the general linear group. For the proof he
could build upon methods developed by I. Schur in his
dissertation from 1901, complemented by an idea of Hurwitz (the
so-called  unitarian restriction) to prove  complete
reducibility. All the irreducible representations could then be
characterized by some {\em symmetry condition} inside some tensor
power $\otimes^k V$. Thus an intriguing correspondence between
the  representations of the symmetric group ${\cal S}_k$
and the irreducible representations of $GL_n (\R)$ inside   $\otimes^k V$
(representations of ``order $k$'')  played an important role in the answer to his first question.\footnote{\cite[455ff.]{Hawkins:LieGroups}}   During the
next two years, this correspondence turned out to be intimately
related to the construction of state spaces for $k$
``indistinguishable particles'' (often electrons bound in an atom) from
the state spaces of the single particles.

This result appeared all the more important to Weyl, as already
{\em before} the advent of quantum mechanics  he had formed  the
conviction that  exactly such irreducible subspaces of $\otimes^n
V$ form the proper mathematical domain of the classical physical
field quantities. He considered the relativistic electromagnetic field tensor
$F^i_j$ with its   antisymmetry property ($n=2$),
\[F^i_j + F^j_i = 0  \; ,\]
as an outstanding example  for this principle.
The methods developed in the study of the general linear group
became the clue to his general theory of representation of the
classical groups.

The {\em second question} had been answered by Weyl already a
little earlier in his investigations of the ``mathematical
analysis of the problem of space''. It  had given him reason to
absorb more of E. Cartan's classification of the infinitesimal Lie
groups than before.\footnote{See \cite{Hawkins:LieGroups},
\cite{Scholz:Infgeo,Scholz:SpaceProblem}. The order of the
questions is  here given according to their relative importance
identified by Tom Hawkins for Weyl's turn towards the new
research project in representation theory of Lie groups.}

During the crucial years 1925 and 1926, Weyl was  busy in other
fields. Immediately after he had finished his  researches in
representation theory of Lie groups, he started intense reading
 for  a book-length article  on    philosophy of mathematics and 
natural sciences, which he had promised to the editors of a
handbook of philosophy.\footnote{Published as \cite{Weyl:PMN}.}
In winter semester 1926/27
he lectured on the theory of continous groups and their
representation as a guest at the  G\"ottingen mathematical institute.\footnote{In this lecture 
Weyl did not yet touch the application of group theory to quantum mechanics \cite{Weyl:VorlesungWS1926/27}. 
I thank  M. Schneider who found H. Grell's  {\em Ausarbeitung} of Weyl's guest lecture in the {\em Nachlass} Herglotz.}
Nevertheless  he was well aware  what was going on in    quantum
mechanics. Even more than that, he actively participated in the
internal discourse  of the protagonists. He was in regular
 communication with E. Schr\"odinger who
taught at the university of Z\"urich in direct  neighbourhood
to the  ETH where Weyl was teaching. And he continued to be
a kind of external ``corresponding  member'' of the G\"ottingen
mathematical science milieu -- notwithstanding his differences
with D. Hilbert on the foundations of mathematics.

\subsubsection*{Communication with M. Born and P. Jordan}
  In  the fall of 1925,  Weyl corresponded with M. Born and P. Jordan on  their actual
 progress  in clarifying Heisenberg's  idea of 
non-commuting ``physical quantities'' in quantum mechanics, which
was initially stated in a mathematically rather  incomprehensible
form.\footnote{\cite{Heisenberg:Umdeutung} submitted July 29,
1925.}
 Heisenberg's idea  was
 ingenious  and opened new  perspectives for theoretical  physics,
but   it was  very difficult to understand. It became  a comprehensible piece of mathematical
physics  only after  the clarification brought about by joint work with Born
and  Jordan   on the one side and by Dirac's contributions  on the other.\footnote{The first paper of Born and Jordan \cite{Born/Jordan}  was   received on September 27, 1925, by  the {\em
Physikalische Zeitschrift}   and a
succeeding one by all the three   
\cite{Born/Heisenberg/Jordan} on November 16, 1925. Dirac joined on November 5, (date of reception) 
\cite{Dirac:qu_m}.}

  Weyl was well informed about the work done by the G\"ottingen physicists and even contributed actively   to  the research discussion among Born,
Jordan, and Heisenberg in the crucial months of mid and late 1925.
In September 1925 Born visited Weyl at Z\"urich and reported him about 
 the latest progress in quantum mechanics. Weyl immediately
started to ``calculate a bit to clarify  things'' for himself, as
he wrote  to Jordan a little
later.\footnote{\cite{WeylanJordan:13Nov1925}} He informed Born
about his insights with great admiration for the work of the
G\"ottingen physicists:
\begin{quote}
Dear Herr Born!\\
Your Ansatz for quantum theory has impressed me tremendously. I
have figured out the mathematical side of it for myself, perhaps
it may be useful for your further progress . \ldots
\cite{WeylanBorn:27Sep1925}
\end{quote}
Weyl proposed to consider the relationship between unitary one-parameter
groups $P(\delta )$ and $ Q(\epsilon )$ with their anti-hermitean infinitesimal
generators $p$, and $q$
 \[ P(\delta ) = 1 + \delta  p +\ldots \;\; \mbox{ and} \;\; Q(\epsilon ) = 1 + \epsilon  q + \ldots  \;\;\;  (0 \leq \delta , \epsilon) \, . \]
He argued that the properties of the (Lie) algebra generated by
pairs of conjugate infinitesimal operators,
\[  pq - qp = \hbar 1, \]
with $1$ the identity and ``$\hbar $ a number'', as Weyl wrote
(he omitted the imaginary factor $i$), could be related to a
commutation relation among the integral operators like
\[  PQ = \alpha  Q P, \;\; \alpha = 1 + \hbar \delta \epsilon  + \ldots \, .   \]
Typical relations among the infinitesimal operators could then be
derived from this approach.\footnote{As an example Weyl presented
the characterization of the formal derivative $f_q:= n p^m
q^{n-1}$ of a monomial $f = p^m q^n$ used by Born and Jordan:
$f_q = \hbar ^{-1}(pf - fp)$.}

 About a week after the submission of his joint article with Jordan,
 Born gave a friendly answer, but with a certain reserve. He   wrote:
\begin{quote}
It was a great pleasure for me to see that our new quantum
mechanics attracts your interest. In the meantime, we have made
considerable progress and are now sure that our approach
covers the most important aspects of the atomic structure. It is
very fine (sehr sch\"on)  that you have thought about our
formulas; we have derived these formulas in our way, even if not
as elegant as you, and intend to publish the subject 
in this form, because your method is difficult  for physicists to
access. \ldots \cite{BornanWeyl:Hs91:488}
 \end{quote}

The communication went on. Weyl received a page proof  of the
submitted paper directly from the  {\em Zeitschrift f\"ur Physik}
and  wrote a supportive letter to the younger colleague, P.
Jordan, in which he apparently referred  to his
alternative approach to the commutation relations once more.\footnote{ On
November 25, 1925, Jordan wrote to Weyl that the latter could
``of course keep the proofs''. In a footnote he  added  an
excuse:   ``I do not know, why they [the page proofs, E.S.] have
been sent to you in such a complicated and demanding form
(umst\"andlich und anspruchsvoller Form). Born and I are innocent
of that  (sind unschuldig daran).'' \cite{JordananWeyl:Hs91:626}.
We can guess that the printer of the {\em Zeitschrift } had sent
the proofs against acknowledgement of receipt, and that  Weyl
was  a bit perplexed by this procedure wondering, perhaps,
whether his G\"ottingen colleagues wanted to make sure their
(undisputed) priority.}

Jordan   thanked Weyl for his comments on November 25, 1925,
shortly after submission of the second paper jointly written with
Heisenberg. He remarked that he had read Weyl's  letter to 
Born  at the time ``with great interest''. He emphasized that  Born
and he had come close to a {\em derivation} of the canonical
commutation relation from the definition of the derivative
$\frac{d}{dt} A$ of an  operator valued function $A = A(t)$ of a
real variable $t$. In a footnote he added:
\begin{quote}
When Born talked to you, we still believed that $pq - qp =
\frac{h}{2 \pi i} 1$ is an independent {\em assumption}.
(emphasis in original)
\end{quote}

  Already in this early correspondence with his  colleagues, Weyl
 looked for unitary groups lying at the  base of the quantization
procedures used by Heisenberg, Born and Jordan. His proposal of
his letter to Born was apparently a first step into
 the direction of using  unitary one-parameter groups obeying a weakened commutativity relation 
(see below, equ.
(\ref{W-commutation})) as a  a  clue to {\em derive} the
Heisenberg relations from basic properties of
projective unitary representations.

In two postcards to Jordan, written in late November 1925,  Weyl
indicated how in his approach an observable $H = H(p,q)$ given in
terms of the conjugate observables $p$ and
 $q$ could be  characterized.
\begin{quote}
I  arrive at   a characterization of the domain of reasonable functions $H$ by the Ansatz
\[ \int \int e^{\xi p + \eta q} \varphi(\xi ,\eta ) d\xi  d\eta\, , \]
which is less formal  than $\sum p^m q^n$.
\cite{WeylanJordan:23Nov1925}
\end{quote}

This was the first indication of what in his publication two years
later \cite{Weyl:QMG} became the proposal to use inverse Fourier
transforms for quantization, the now so-called Weyl-quantization
(equations (\ref{Fourier}) and (\ref{W-quantization}) below).
 Born and his assistant Jordan decided, however, that Weyl's
 approach was too cumbersome for the introduction of the new quantum
 mechanics to the physics community, and relied on their own approach.
 The long delayed and selective  reception of Weyl's idea  shows
that Born may have been  right in this estimation. On the other hand,   his decision may have contributed
to  the long delay for a recognition of Weyl-quantization as a useful approach 
in mathematical physics.

\subsubsection*{Abelian ray representations}
Weyl came back to his early proposals nearly two years later in
his first article dealing with quantum mechanics
\cite{Weyl:QMG}.\footnote{Received October 13, 1927.}
 He clearly distinguished between {\em pure states} and of {\em mixtures}. Pure states were  mathematically represented by eigenvectors (or more precisely by corresponding   complex  unit rays) of the typical observables which described the defining properties of a particle or dynamical state. Mixtures, on the other hand, were described   contextually as
composed from pure states in ``any mixing ratio''
\cite[97]{Weyl:QMG}.  In this way Weyl indicated that a mixed
state might be characterized by a probability measure on the
state space,  although he did not spell  out details.
A little later, and originally  without knowledge of  Weyl's
manuscript, von Neumann proposed to  formalize both    mixed and
pure systems  by (positive) hermitian operators $A$. Pure states
were  those given by projection operators onto one-dimensional
subspaces and  mixtures by more general positive hermitian
operators  \cite[215ff.]{Neumann:WAQ}.\footnote{Von Neumann
presented his paper on November 11, 1927, to the {\em G\"ottinger
Gesellschaft}. In the page proofs he added a reference to Weyl's
paper  \cite[219, footnote]{Neumann:WAQ} and vice versa \cite[90,
footnote]{Weyl:QMG}; compare \cite[431ff.]{Mehra/Rechenberg:VI}.
In later terms, von Neumann's positive hermitian operator $A$
can  be related to a {\em trace class} operator $T$ by $A =
(T^{\ast} T)^{\frac{1}{2}}$, where $T$ is  of unit trace norm $T_1= 1$. Here   $|T|_1:= Tr T = \sum_{k} (T u_k, u_k) = 1$ with respect to any  complete
orthonormal set  $\{u_k \}$. Moreover, the trace of $T$ can be
calculated by the sum of the (positive) eigenvalues $a_{\nu}$ of
$A$, \quad $Tr T = \sum_{\nu} a_{\nu}$.}

Weyl's main point was, however, the discussion of what he
considered the ``more profound''  question of the ``essence
(Wesen) and the correct definition of canonical
variables''\cite[91]{Weyl:QMG}  $P$ and $Q$, satisfying the {\em
canonical} or {\em Heisenberg commutation} relation: \beq
\label{H-commutation} [P, Q] = \frac{\hbar}{i} 1.  \eeq He
proposed to relate any hermitian operator $A$ to  the unitary
1-parameter group generated by its skew hermitian relative $iA$
\[  t \mapsto e^{i t A}  \]
and to consider the quantum mechanical observables from an
``integral'' point of view, in the sense of the generated
1-parameter groups. That was a conceptual move  similar to the one in Weyl's work on
representation theory, where he found intriguing new aspects by
passing from the infinitesimal point of view (the Lie-algebras in
later terminology) to the integral perspective (the groups
themselves).

Turning the perspective round, he considered a classical state
space described by pairs of $n$ conjugate observable quantities
$(p,q)$, such as the spatial displacement $q$ with respect to a
frame and its conjugate  momentum  $p$. Then the  state  space
could be considered as
 an abelian group $G$ of two continuous parameters $(t,s) \in   \R^2 = G$ (in
 the case of $n=1$ pairs). For the quantization it was natural to
look at a {\em  unitary  ray representation}, i.e. a
representation up to multiplication by a complex number of unit
norm.

Then it was clear that in the quantum context  the  commutation relation for the
generating  1-parameter groups $e^{it P}$ and $e^{is Q}$  have to be weakened. Commutativity had to hold only up to a unitary factor, 
\beq \label{W-commutation} e^{i s P} e^{ i t Q} = e^{ic \, st }
e^{i tQ} e^{i s P} , \eeq 
where $c$ is a real constant normalized to
$c=1$ or $c= \hbar$. Let us refer to equation
(\ref{W-commutation}) as the {\em Weyl-commutation} relation for
conjugate pairs of 1-parameter groups in unitary projective
(quantum) representations.

Weyl showed that for the corresponding  skew-hermitian infinitesimal generators $iP$, $iQ$ the  deviation
(\ref{W-commutation}) from strict commutativity   implies
\[  PQ - QP = - i c \, 1 , \]
i.e., the Heisenberg commutation rule (\ref{H-commutation}) for  a pair of conjugate observables.

Weyl generalized this procedure to $n$-tuples of pairs of
observables $P_1, Q_1,$ $ \ldots  ,P_n,Q_n$. Then a representation
on quantum rays\footnote{``Quantum ray'' signifies that from
the one-dimensional subspace, the classical projective ray, only
the norm 1 representatives play a role in the quantum
mechanical context. }
 allowed to modify  the strict commutation relation of an
 abelian group $(t_1, \ldots, t_n, s_1, \ldots s_n) \in G =
\R^{2n}$ to slightly deformed  Weyl-commutation relations of
the form
\[ e^{is_{\mu} P_{\mu}} e^{it_{\nu} Q_{\nu}} = e^{ic \delta ^{\mu}_{\nu}\, s_{\mu} t_{\nu} } e^{it_{\nu}Q_{\nu}} e^{is_{\mu} P_{\mu}} , \]
with $\delta ^{\mu}_{\nu}$ the Kronecker delta  and $c=1$, or $c=
\hbar$. For the infinitesimal generators that corresponded to a
normalized form of the skew symmetric system of coefficients
$c_{\mu \nu}$ in  the system of relations \cite[114]{Weyl:QMG}
\beq \label{W-comm-inf}  P_{\mu} Q_{\nu} - Q_{\nu} P_{\mu} = - i
c_{\mu \nu } \, 1 .\eeq 
That  led to intriguing relations for the
addition rule for the $2n$-parameter  unitary ray representation.
If we use the  denotation  $(s,t) \in \R^{2n}$  and
\[  W_{s,t}:= e^{is_1 P_1} e^{is_2 P_2} \ldots e^{is_n P_n} e^{it_1 Q_1}  \ldots e^{it_n Q_n}  ,\]
the addition becomes
\[   W_{s+s', t+t'} =    e^{- i c<s',t>} W_{s,t} W_{s',t'} \; , \]
where $<s',t> := \sum_{\nu} s'_{\nu} t_{\nu}$ and, as above,
$c=1$ or $c= \hbar$. The resulting structure was an {\em  irreducible
projective unitary representation of the abelian group} $G =
\R^{2n}$;   Weyl  called it  an  ``irreducible abelian  rotation
group operating on a the field of rays (Strahlenk\"orper) of pure
states '' \cite[118]{Weyl:QMG}. He restricted his   investigation
to the case of everywhere defined, bounded (skew-) hermitian
generators  and the resulting unitary transformations and gave a
sketchy argument that these were the only irreducible projective
representations for each $n$.

For a serious application to quantum mechanics, the
generalization to the  case of unbounded operators was, of
course,  important. It   was solved independently   by  Marshall
Stone and von Neumann \cite{Stone:Unitary,Neumann:UnitaereOp}.
Von Neumann showed,  in addition, that the Weyl commutation
relations ((\ref{W-commutation}), (\ref{W-comm-inf})) characterize
irreducible  unitary projective representations of continuous
abelian groups up to unitary isomorphism.

\subsubsection*{Weyl quantization}
Weyl, on the other hand, continued his article by looking for  a   procedure which could give
 operator companions to  (classical) physical quantities in a systematic way, i.e., he looked for a systematic approach to {\em quantization} 
 \cite[116]{Weyl:QMG}. If a classical quantity is expressed  by  a
function $f(p,q)$   of the canonical variables $p, q$ ($f \in
{\cal L}^2 \R^{2}$ for $n=1$), he looked at  the  Fourier
transform $\xi$ of $f$. Then $f$ can be gained back from $\xi$ by
\beq \label{Fourier}  f(p,q) =   \int  e^{i (ps + qt)}  \xi (s,t)
ds dt  , .\eeq
 Weyl proposed to use the analogously formed operator-valued integral
\beq \label{W-quantization} F:=  \int   e^{i (Ps + Qt)} \xi
(s,t)  ds dt  = \int \xi \,W_{s,t} \, ds dt \eeq as the quantum
mechanical  version of the physical quantity related to $f$. In
case of periodic variables,   pairs $(p,q)$  represent elements
on the  torus  $ G = T^2 : = S^1 \times S^1 \cong \R^2/ \Gamma$, where
$\Gamma$ is the lattice generated by the periods. Then  the
integration reduces to a summation over integer numbers $s$ and
$t$ in $\Z$, because the Fourier transform $\xi$ lives on the
discrete domain $\hat{G} = \Z^2$.  Moreover,  $f$  is an element
of the function algebra on the abelian group $G= \R^2$, or $T^2$
in case of periodic variables. For a real valued function $f$, in
particular, the corresponding $\xi$ satisfies
\[  \xi (-s,-t) = \overline{\xi} (s,t) \]
and leads to a hermitian operator $F$.

 In the methods introduced and used by physicists at the time for the
 quantization of classical observables, $ p \mapsto P$, $q \mapsto
Q$ , the non-commutativity of $P$ and $Q$ led to a fundamental
difficulty for  an observable given as a  function $f(p,q)$ of
the basic dynamical variables $p$ and $q$. Already in the simple
cases of a polynomial function, it was not clear which operator
one should choose for the formal expression $f(P,Q)$. For example  for
$f(p,q) = p^2q$ one could choose any of  $P^2Q$,  $PQP$ or
$QP^2$, etc.. Weyl's  unitary ray representation approach
resolved (or avoided) this difficulty from the outset. The
operator inverse of the Fourier transform (\ref{W-quantization})
gave a unique and structurally well determined assignment $f
\mapsto F$  of hermitian operators to real valued quantities.
Weyl was  therefore convinced   that  ``our group theoretic
approach  shows immediately the right way'' towards the
quantization problem \cite[117f.]{Weyl:QMG}.

Of  course, the whole approach worked only for non-relativistic
mechanical systems in which time is ``the only independent
variable'', whereas field theory deals with quantities extended
over time {\em and} space, which relate observations and
measurements among each other.
 Weyl considered the independent variables as ``projected into the world''
 by arbitrary conventions in such a manner that the dependence of
physical quantities on them could {\em not} be measured
\cite[124]{Weyl:QMG}.  In this sense,   the independent variables
played for him the role of some kind of a-priori component in
theory construction. They  were necessary for the conceptual
architecture of the whole symbolic construction, although they
were not directly related to observable quantities. In
non-relativistic quantum mechanics time was the only
``independent variable''  left. He added:
\begin{quote}
If one wants to resolve the criticized omission of the time
concept of the old pre-relativistic mechanics,  the observable
quantities time $t$ and energy $E$ have to be considered as
another canonically conjugate pair, as is indicated already by
the action principle of classical mechanics.  The dynamical law
[of the Schr\"odinger equation, E.S] will then completely
disappear. \cite[127]{Weyl:QMG}
\end{quote}
He referred to Schr\"odinger's first attempts   to obtain a
relativistic theory of the electron in a centrally symmetric
field,
 but neither here, nor in any later publications, did he
start to work out this  idea of how one might proceed to
build  a relativistic quantum field theory. A
good occasion would have been his contributions to Dirac's electron
theory, two years later; but by then  he had already accepted that
the physicists working on this question -- Dirac, Jordan,
Heisenberg, and Pauli -- had chosen a completely different approach.
They developed the method of   so-called second quantization,
which seemed  easier ``to access for physicists'', to take up  Born's words from   his letter of October 2, 1925 to Weyl. 

The problems sketched  in Weyl's 1927 paper,  the method of
unitary  ray representations of commutative groups, and the
ensuing quantization method proposed were soon reconsidered in
Weyl's book \cite{Weyl:GQM} and made more accessible to an
international audience by its English translation in 1931. The only traces it left on contemporary work was that of 
von  Neumann and  Stone, mentioned above. But it turned
out to be of long range inspiration. In the next
generation, G. Mackey took up Weyl's representation theoretic
perspective  and  developed it into a broader program for the study of
{\em  irreducible projective representations} as a starting point
for a more  structural understanding of quantum physical systems
\cite{Mackey:1949}.

In the 1960s, Weyl's quantization started to be revitalized. In
this decade, the torus case, $G = T^2$, was reconsidered as a
special, and the historically earliest, way to  introduce  a
deformed product on the  Fourier dual group, $\hat{G} = \Z^2$.
For two elements $f, h$ of the function algebra on $G$ with
Fourier transforms $\xi = \hat{f}, \; \eta = \hat{h}$,  $\xi,
\eta \in \hat{G}$, let the Weyl quantization be written as  $f
\mapsto F$, $h \mapsto H$. Then the  composition  of the Weyl
quantized operators
\[F \cdot  H \]
 could be transported back to the  original functions $f, h$ or their Fourier
 transforms $\xi, \eta$. That  led to  a {\em deformed product}
depending on a parameter $c$ (typically $c = 1$ or $c = \hbar$),
\[ f \ast_{c} g , \;\;\; \mbox{respectively} \;\; \xi \hat{\ast }_{c} \eta  \, ,\]
with properties which attracted a new generation of
researchers.\footnote{For an overview see
\cite{Rieffel:Quantization}. Rieffel refers to
\cite{Pool:WeylCorr} as the first paper in which an explicit
description of the deformed product on the Fourier transform
functions was given.  His claim that already von Neumann had
``pointed out that Weyl quantization induces a new
product on functions'' \cite[70]{Rieffel:Quantization}  seems, however, to be anachronistic. The closest approximation
to such a view in von Neumann's paper is, as far as I can see, a
reference to the ``Gruppenzahlen''  at the end of the paper,
where the  terminology ``Gruppenzahlen'' refers to  functions $f$
on $G$  as elements of the group algebra $\C[G]$
\cite[229]{Neumann:UnitaereOp}. Such a perspective  was also
discussed  in Weyl's paper \cite[106]{Weyl:QMG} (and there even in
more  detail). In the abelian case considered here the group
algebra is commutative and could at best serve  as the starting
point for the introduction of the deformed product.  Neither von
Neumann nor Weyl  mentioned the idea that the Weyl-quantized
operators might be used to introduce a modified (non-commutative)
product of the ``Gruppenzahlen'' themselves.}

The resulting non-commutative function algebra on the torus $T^2$
or its Fourier dual $\hat{T}^2 = \Z^2$ became the starting point
for the study of the non-commutative torus, one of the first
well-known  cases of non-commutative geometry.  Weyl-quantization
turned out to be just one among  a larger class of deformation
quantization procedures.

Thus Weyl's first paper presented  ideas to the public, which  he
had developed essentially when he was  still  ``at the backstage''
of the  quantum mechanical scene, as we have called it, turned
out to have long range impact  in several respects,
\begin{itemize}
\item[---] for the study of {\em  irreducible projective representations} (Stone, von Neumann, Mackey e.a.),
\item[---]   as an inspiration for the search for conceptually founded
quantization procedures such as the {\em Weyl-quantization}, as it
was called after the 1960s,
\item[---]  and finally as one of the sources for  a  {\em non-commutative modification} of the
 the torus (Pool, Rieffel e.a.).
\end{itemize}
At the time of their publication, Weyl's proposals were,
however,  far too distant from contemporary quantum mechanical
research to be taken up in the physics community.
For several decades the paper \cite{Weyl:QMG} remained a lonely standing  monument.

\section*{5. Weyl entering the stage }
In late 1927, Weyl entered the field of quantum mechanics with
full force. He had announced a lecture course on group theory at
the Z\"urich  {\em Eidgen\"ossische Technische Hochschule}, ETH,  for winter semester 1927/28.
In the  summer of this year,
both Z\"urich theoretical physicists   accepted calls to other
places, E. Schr\"odinger left the University of Z\"urich and went
to Berlin; P. Debye gave up his chair at the  ETH on occasion
of a call  to Leipzig. Weyl used the opportunity to reorient his
lecture course originally  announced on group theory only and offered
it now as a course on  ``Group theory and quantum mechanics
(Gruppentheorie und Quantenmechanik)", without running the risk of
putting off his local colleagues in physics. Now he had a good
opportunity to present his views on group theoretical methods in
quantum mechanics. His main interest was centered on  the intriguing interplay between
representations of the orthogonal group $SO_3$ (and $SU_2$) and
the permutation group, which  about the same time Wigner and von Neumann  hit upon from their side.
 Let us remember that in
summer or autumn 1927 only Wigner's own papers were published. The
joint work with von Neumann was still going, on when Weyl prepared
the book manuscript from the lecture notes in the summer semester
1928. In late August the book was finished and given to the
publisher. In the sequel we will also use the abbreviation GQM for
it \cite{Weyl:GQM}.\footnote{If not otherwise stated, quotations
refer to the first edition of GQM. If possible translations are
taken from H.P. Robertson's English version of the second
edition; where necessary or advisable (because of meaning affecting
shifts) direct translations from the first edition are given by
the author (E.S.). The second edition will be quoted by
\cite[$^2$1931]{Weyl:GQM}, the English translation by
\cite{Weyl:GQMEnglish}. For a discussion of the book see
\cite{Speiser:Kiel}.}

Weyl's contributions  to the topic and the joint work by Wigner and von Neumann were developed
in parallel and independently of each other,  as far as any direct
exchange of ideas is concerned. They nevertheless established a
common theoretical approach to groups in the quantum mechanical
explanation of atomic spectra. This is a good  case for a
comparative study of   how Weyl's
perspectives as a mathematician with great expertise in group
representations influenced his approach to the subject. We can
compare it directly with the Wigner -- von Neumann ``team'', one
of them  (von Neumann) a brilliant mathematician who had
assimilated the new results in representation theory in a speed
which later became legendary, the other one  a theoretical
physicist of admirable mathematical powers.

Two  points of the broader story of group
theoretical methods in quantum physics have to be mentioned, before we come to the discussion
of Weyl's treatment of the interplay of the symmetric and the
orthogonal groups in spectroscopy and quantum chemistry. Here we can
only mentioned them in passing, although they  deserve closer
scrutiny in their own contexts. 

\subsubsection*{General relativistic spinor fields}
 Exactly at the end of Weyl's course and shortly after it
 finished, Dirac's two path-breaking  papers on the relativistic
theory of the electron appeared  \cite{Dirac:1928_I/II} and found
immediate recognition \cite{Kragh:Dirac}. Therefore Weyl's book
already contained a chapter on Dirac's theory. Later in the year
1928 and early the next one, Weyl  took up Dirac's theory,
simplified it from the point of view of group representations and
put  it into a   general relativistic framework. For physical
reasons, Dirac worked with a reducible representation of the
Lorentz group, now written as $D^{(\frac{1}{2},\frac{1}{2})}$,
whereas Weyl proposed a reduction to irreducible components,
characterized by the standard representation of $SL_2\C$ in
$\C^2$,  $D^{(\frac{1}{2},0)}$, and/or its conjugate
$D^{(0,\frac{1}{2})}$ (``Weyl spinors'' versus ``Dirac spinors'',
in later terminology). Weyl's main goal in a series of papers in
the year 1929 was,  of course, of a different nature, the
adaptation of spinor theory to general relativity. In this
enterprise he had again independent parallel workers, V. Fock and
D. Ivanenko at Leningrad.  Weyl and Fock/Ivanenko built
essentially the same core theory, but differed in outlook and
details. That is an interesting story in itself, which cannot be
told here.\footnote{Compare
\cite{Vizgin:UFT,Goenner:UFT,Straumann:DMV,Scholz:Fock/Weyl}.}
Weyl did {\em not  include} this generalized treatment of the
Dirac equation in the second edition of the book, but only
referred to it in passing  at various places \cite[$^2$1931, VII,
195]{Weyl:GQM}.

In the second edition  he changed and  extended, the treatment of the
 special relativistic Dirac equation. In the
first edition he  discussed a non-relativistic first approach
to  ``second'' quantization  of the electron and the electromagnetic
field  \cite[\S 44]{Weyl:GQM}. At the end of the passage Weyl
remarked:
\begin{quote}
We have thus discovered the correct way to quantize the field
equations defining electron waves and matter waves. The exact
realization will be the next task of quantum physics; the
maintainance of relativistic invariance seems to offer serious
difficulties. Here again  we find that quantum kinematics is not
to be restricted by the assumption of Heisenberg's specialized
commutation rules. And again it is group theory, which supplies
the naturally generalized variant, as is shown by the next
section \ldots [in which unitary ray representations and the first
steps of  Weyl-quantization  were presented, E.S.]. 
\cite{Weyl:GQM} \cite[$^2$1931, 203]{Weyl:GQM}
\end{quote}

In summer 1928, he apparently still assumed that his approach to quantization might allow a generalization from the group $\R^3$ of non-relativistic  kinematics to the relativistic case.
 In the second
edition he omitted the second and the last sentences, after   in  January 1929
 Heisenberg and Pauli had made decisive progress in their approach to ``second
quantization''. Weyl still kept the passage on unitary ray
representations and  to (Weyl-) quantization, but  no longer recommended 
  his own approach as a a path towards  relativistic field
quantization. He included two new sections with  a discussion
of this new and difficult terrain, following Pauli, Heisenberg
and Jordan, although now the obstacle of uncontrollable inifinities appeared at the horizon
 \cite[$^2$1931, chap.IV, \S \S 12, 13]{Weyl:GQM}. 

\subsubsection*{Discrete symmetries}
In
these new passages Weyl started also  to explore  the role of discrete symmetries in the context of early
relativistic field theory,  parity change $P$,  time
inversion$T$, and  charge conjugation $C$. They ended  with a remark   which
struck readers of the next generation as  surprising 
 and even ``prophetic'':
\begin{quote}\ldots this means that positive and negative electricity
 have essentially the same properties in the sense that the laws
governing them are invariant under a certain substitution which
interchanges the quantum numbers of the electrons with those of
the protons [later readers would functionally rephrase the term
by ``positrons'', E.S. ]. The dissimilarity of the two kinds of
electricity thus seems to hide a secret of Nature which lies yet
deeper than the  dissimilarity of past and future.
\cite[$^2$1931, English, 264]{Weyl:GQM}
\end{quote}
We cannot take up the thread of the rise and establishment of the
discrete symmetries in quantum field theory here;  readers interested in this topic  may
like to have a look at the discussion in
\cite[293]{Coleman/Korte:DMV} and \cite[141]{Straumann:DMV}.

\section*{6. Weyl on stage }
We come back to comparing the 
 differenct  outlooks of Weyl and   Wigner/von Neumann on  groups in quantum mechanics. Technically, they agreed completely, as Weyl frankly stated when he wrote the preface to
his book in August 1928.\footnote{Remember that all three parts of the Wigner/von
Neumann series had appeared at that time, the last one in June 1928.}
Discussing the role of group representations in quantum
mechanics, he observed:
\begin{quote}
The course of events is so inevitable (zwangsl\"aufig) that nearly
everything that was still new at the time when I gave the course
has been published elsewhere in the meantime, in particular by
the work of the colleagues (der Herren) C.G. Darwin, F. London,
J. von Neumann and E. Wigner.
\end{quote}
He added:
\begin{quote}
That is different with Dirac's wave equation of the electron,
which introduced essential new ideas into the theory during the
time when this book was being written. \cite[vi]{Weyl:GQM}
\end{quote}
The reference to F. London, and at other places to W. Heitler, referred to the
theory of molecular bonds, which
Weyl had approached with the tool kit of representations of the
symmetric group, starting from the joint article of Heitler and
London.\footnote{\cite[$^2$1931, 300, chap. V, endnote 10]{Weyl:GQM}.  Darwin's publications dealt with the spin phenomenon; among them \cite{Darwin:electron_I,Darwin:electron_II}. It did not involve explicit group theoretic aspects. }
Even more than the other authors, Weyl emphasized the structural role group
representations  for the understanding of
quantum physics. He hoped that they would survive future changes  of the actual mathematical models of the
atomic or molecular systems:
\begin{quote}
Recently it turned out that group theory is of fundamental
importance for quantum mechanics. In this context it reveals the
most essential features {\em whatever the form of the dynamical law may be, i.e., without definite assumptions on the forces which are acting.}
\cite[2, emphasis E.S.]{Weyl:GQM}
\end{quote}

The last remark described  
quantum mechanics as a theory in development. Weyl considered  it to be
in an   unfinished  state. That differed from  the
credo of the Copenhagen -- G\"ottingen group which argued
strongly in favour of having  achieved a ``completion'' of
quantum mechanics.\footnote{Compare the  title of volume VI of \cite{Mehra/Rechenberg}: ``The Completion of Quantum Mechanics 1926 -- 1941''.}  Weyl did {\em not} share, however,
Einstein's opinion that quantum mechanics had to be considered as of only
provisional character, as long as its purely stochastic
determination  was not reduced to a classical field theory lying at its base.  Weyl even  had welcomed
the stochastical character of natural laws well before the birht of the ``new''
quantum mechanics \cite{Weyl:Kaus_Stat}.  Of course, he was,
well aware of the fundamental  problem that quantum
mechanics and relativity had established two theories of  basic
levels of nature, which were 
conceptually and mathematically far apart. Already during his   ``backstage
period'' Weyl had  looked for possibilities of reconciliation of
relativity theory and quantum physics  (see above). In
summer 1928, after Dirac's
breakthrough to a first relativistic quantum theory with
empirical successes,  he expected further changes to come. In such a period, Weyl thought that   the assumptions on the ``form of the
dynamical law'' might still be subject to considerable  change. The representation theoretical methods, on the other hand, appeared to him as part of a stable core of quantum mechanical knowledge.

This  conviction of a deep structural meaning of group
representations was the central topic in  GQM. Similar to his
first book on mathematical physics, {\em Space - Time - Matter},
Weyl gave a complete introduction to the mathematics of the field
and wrote one of the first textbook expositions of quantum mechanics. He
started with an introduction to what he called  {\em unitary
geometry}, i.e., the theory of  Hilbert spaces and the
diagonalization of hermitian forms,  although  essentially
restricted to the finite dimensional case (chapter I). He
continued with  an introduction to quantum mechanics integrating
the Schr\"odinger view of the dynamical law in the
non-relativistic case  and the G\"ottingen
(Heisenberg-Born-Jordan) point of view of observables represented
by hermitian operators and their  quantum stochastical
interpretation (chapter II). Of course, he emphasized the turn
quantum mechanics had  taken with respect to  classical natural
science. Both had in common to be ``constructive''.
\begin{quote}
Natural science is of a constructive character. The concepts with
which it deals are not qualities or attributes which can be
obtained from the objective world by direct cognition. They can
only  be determined  by an indirect methodology, by observing
their reaction with other bodies; their implicit definition is
consequently conditioned by definite laws of nature governing reactions. \cite[66]{Weyl:GQM}
\end{quote}
Classical mechanics was able to  assume that
such ``constructive properties'' were attributes of the ``things
as such  (Dingen an sich)'', in the sense of pertaining to them,
even if the manipulations necessary to their determination were
not undertaken. In    quantum 
physics this  was no longer possible. In this point   Weyl  agreed with N. Bohr.
\begin{quote}
{\em With  quanta we run into a fundamental barrier (Schranke) to
this epistemological position of constructive natural science.}
(ibid., emphasis in original, my translation, ES)
\end{quote}
This  limitation lay at the basis of Heisenberg's undeterminacy
relation.  Weyl accepted it as a fundamental insight, different from Heisenberg's  mathematical characterization of the commutation relation.\footnote{Weyl presented  Heisenberg's undeterminacy
  in a form due to a  communication
by W. Pauli \cite[67, appendix 1]{Weyl:GQM}.}

In the third section Weyl introduced  the representation theory
of finite groups with some general remarks on continuous groups,
their characters and their infinitesimal groups (chapter III).
The presentation of concrete examples, in particular the
orthogonal group, Lorentz group, the special unitary  and the
symmetric groups were postponed  to the later sections on
``applications of group theory to quantum mechanics'' (chapters
IV and V). Chapter IV   contained the theory of atomic spectra, Dirac's electron theory, and
his own method of unitary ray representations. The last chapter
developed the combined theory of representations of the unitary
group and the symmetric group, preparing his approach to  the
theory of valence bonds (chap. V).

His presentation of atomic spectra \cite[157ff.]{Weyl:GQM} relied  much more on theoretical  arguments  and used less explicit calculations of
eigenfunctions than  Wigner/von
Neumann's. Nevertheless his discussion went  as deep  into the
physics  context as Wigner's. It included, among others,  a concise group
theoretic discussion of Pauli's  mathematization of spin and of  the  
anomalous Zeeman effect. Weyl apparently wanted to demonstrate  the usefulness of the  structural view of mathematics  for a conceptual understanding in  physics.

\subsubsection*{Pauli spinors from the point of view of representation theory }
For the characterization of electron spin Weyl could
build upon his observation of 1924, that the
special orthogonal groups $SO_n \R$ are not simply connected but possess, for $n>2$, a two-fold universal covering group  \cite{Weyl:62}.  He
 clearly   distinguished ``two-valued'' and
one-valued representations of these groups 
\cite[II, 602ff.]{Weyl:Darstellung}. 
For the introduction of electron spin, he    nevertheless preferred the
 more physical approach of extending Schr\"odinger wave
functions to Pauli spinors. To concentrate ideas, he started with
the discussion of alkali spectra, governed by one external
electron with a state space called ${\cal E}$:
\begin{quote}
We deal with a single electron; the wave function depends only on
$t$ and the three space coordinates $x, y, z$.  It cannot be a
scalar, however, but is a two-component covariant quantity of
type ${\cal D}_{\frac{1}{2}}$. Then we have ${\cal D}= {\cal
D}_{\frac{1}{2}} \times {\cal E}$, and the decomposition of
${\cal E}$ into its irreducible components ${\cal D}_l$  with the
integer azimuthal quantum number $l$ gives the constituents
${\cal D}_{\frac{1}{2}} \times {\cal D}_l$. Each of those
decomposes again into a doublet ${\cal D}_j$ with $j= l +
\frac{1}{2}$ and $j = l - \frac{1}{2}$ \ldots .
\cite[162]{Weyl:GQM}\footnote{Weyl's ${\cal D}_j$ corresponds, of course, to our ${\cal D}^{(j,0)}$ of equation (\ref{representation SU2}).}
\end{quote}

The  observation of the last sentence was an immediate
consequence of the decomposition formula for a tensor product of
representations of  $SU_2$, given here in Weyl's notation
\cite[166]{Weyl:GQM}
\[ {\cal D}_s  \otimes {\cal D}_l  = \sum_{j=|l-s|}^{l+s}{ \cal D}_j  \; .  \]

 As the old theory without spin characterized the terms very well up to small
 effects, Weyl assumed 
 that the
two-component wave functions were  well approximated by   the
``old''  Schr\"odinger wave functions  (as did his quantum physical colleagues). The dimension of
the function space was now doubled, with a corresponding
rise in the degree of degeneracy. He introduced the notation ${\cal E}_l $ for an invariant subspace of ${\cal E} $, ${\cal E}_l \cong {\cal D}_{\frac{1}{2}} \otimes {\cal D}_l    $ and gave his interpretation of the appearance of spin doublets:
\begin{quote}
\ldots thus ${\cal E}_l$ now possesses all pairs $\psi = (\psi_1,
\psi_2)$ as eigenfunctions \ldots. They obviously form a linear
manifold of $2 (2l+1)$ dimensions. But now a small perturbation
term will be added to the wave equation, the
``spin-perturbation'' which couples the components $\psi _1, \psi
_2$ among each other. Thus the former accidental degeneracy is
broken, the $2 (2l+1)$-fold eigenvalue $E_l$ is split into two
values of multiplicities $2 j +1$, with $j = l \pm \frac{1}{2}$,
just as the representation ${\cal D}_{\frac{1}{2}} \times {\cal
D}_l$ is decomposed into two irreducible constitutents. This is
the theory of the doublet phenomenon as sketched by W. Pauli.
(ibid.)
\end{quote}

This was a beautiful demonstration of how representation
theoretic structures  appeared very naturally  in the
material of basic quantum mechanics. They were able to elucidate
the symbolic constructions and the perturbation arguments introduced by contemporary physicists,
including the kind of structural approximation which led from 
Schr\"odinger's  to the Pauli's  wave functions.

In the discussion of the anomalous Zeeman effect,  i.e, the split
of spectral lines of multiplets  under the influence of an
external magnetic field,
 Weyl showed that the representation theoretic view could also lead to
 quantitative results; he gave a theoretical derivation of the {\em Land\'e formula}
 for the split of  spectral terms  in an external magnetic field
\cite[164ff.]{Weyl:GQM}.\footnote{Land\'e had determined a characteristic 
factor $g$, important for the calculation of the widths of   the
line split,  as $g = \frac{2 j +1}{2 l+1}$, where $l$ was the old
(integer valued) azimuthal quantum number and $j = l \pm
\frac{1}{2}$ an ad-hoc  modification which could later be
interpreted as the ``internal'' quantum number of  the
representation ${\cal D}^{(j,0)}$ , taking spin into account. Weyl
derived $g$  in very good approximation  from the magnetic
momenta of the Pauli-spinors as $g-1 = \frac{j(j+1)- l(l+1)+
\frac{3}{4}}{2j(j+1)}$, which reduces to Land'e's formula in the
cases $j= l \pm \frac{1}{2}$. Compare
\cite[499]{Mehra/Rechenberg:VI}.}

\subsubsection*{A physical role for representations of the symmetric group}
In his presentation of molecular bonds and its group theoretic
background (chap. V), Weyl was apparently intrigued
by a structural  analogy of the
 spin-coupling problem of the $n$-electron system with his general studies of group representations.
 In both cases, a strong and deep interplay of a continuous group
($SO_3$ or $SU_2$ in the spin case, more generally
 any classical group)  with the operation of the symmetric group, or some subgroup (the
Weyl-group in the general case), formed the essential core of his
analysis. Thus Weyl declared that one of the goals of his lecture
course and the book was to give a unified picture of the
representation theory of finite and of continuous groups.
\begin{quote}
Already from the purely mathematical point of view, it no longer
seems justified to make such a sharp distinction between finite
and continuous groups as is done in the traditional textbooks.
\cite[V]{Weyl:GQM}
\end{quote}

He was very pleased that the  study of the spin of an
$n$-particle system relied on what he called at different occasions a  {\em bridge}
between the discrete and the continuous group representations \cite{Weyl:87}. His goal was to make
this bridge conceptually  as clear as possible, not only to use
its consequences in the determination of term systems or in the investigation of  chemical bonds.  This does not mean that he contented himself with purely structural insights. He rather started to elaborate
the representation theory of the symmetric group  with the
explicit goal to derive calculatory tools.  For this purpose he  refined the use of Young diagrams and Young tableaus.

In the last respect he made considerable advances after the
publication of the book. Several articles on this
topic followed during the next year, among it the main research paper
\cite{Weyl:79} and some  expository ones
\cite{Weyl:80,Weyl:86,Weyl:87}. In these papers Weyl achieved a  structural clarity in  the study of spin-coupling, comparable to
the one he had gained  during the years  1925/26 for the
representation theory of the classical groups. On the basis of
these results he  completely rewrote the last part of his book
(chapter V) for the second  edition (and its English
translation). The revised chapter V  became the source for a
tradition of a long, although slow, trickling down  of knowledge
and of symbolical tools from the representation theory of the
symmetric group to the  theory of atomic and molecular
spectroscopy (later even to nuclear spectroscopy) and to quantum
chemistry.

In these considerations Weyl employed 
similar methods to those he had developed in his studies of 
representation theory  in 1924/25. Central for both approaches
was the association of a {\em symmetry operator} $A$  to each element  $a$ of the group
algebra  $\C [S_f] $ of the symmetric group ${\cal S}_f$, operating on a tensor
product space $\bigotimes^f V $. Using
Weyl's notation $F= F(k_1, \ldots, k_f) $ for a tensor $F \in
\bigotimes^f V$,\footnote{This notation takes allows to use a shorthand notation for 
 the operations of $\C [S_f] $ on general  tensors $F = \sum_j \alpha_j \, v_1^{(j)} \otimes \ldots \otimes v_f^{(j)}$, defined by linear extension of the naturally defined operation on the decomposable tensors $ v_1^{(j)} \otimes \ldots \otimes v_f^{(j)}$.}
 the symmetry operator $A$ associated to 
\[a = \sum_{s \in S_f} a(s) s \;  \in \C[{\cal S}_f ] \] 
was given by:
\[ A:  F(k_1, \ldots, k_f) \mapsto \sum_{s \in S_f} a(s)  F(k_{s(1)}, \ldots, k_{s(f)})  \; . \]
Using such symmetry operators, Weyl formulated symmetry conditions for elements in the
tensor space   $\bigotimes^f V $ and showed that invariant subspaces of the regular representation on $\C[{\cal S}_f]$ specify invariant subspaces of $GL (V)$ on $\bigotimes^f V $. 

\begin{theorem}
 There is a $1:1$ correspondence between invariant subspaces of the regular
 representation of $S_f$ and invariant subspaces of the operation
of $GL(V)$  on  $\bigotimes^f V $. The same holds for its 
irreducible building blocks (the corresponding irreducible representations). \\
{\em \cite{Weyl:79}, \cite[350]{Weyl:GQMEnglish} }

 \end{theorem}

A comparable  correspondence had already been used by I. Schur in
his dissertation \cite{Schur:Diss} and, in a modified form again
in \cite{Schur:1927}. Weyl gave  full credit to these works. Only his
method of symmetry operators was new, and he thought it 
to be of advantage for the clarification of the overall structure
of the correspondence. In an exchange of letters, which is only
partially preserved,   Schur expressed complete
consent:\begin{quote} I do not find anything in your interesting
paper which I had to object to. I even accept   as not
illegitimate the gentle criticism which you offer to my
publication from the year 1927. I am very glad to see that you
emphasize the connection between my old approach from the year
1901  and your elegant formulation. I also give preference to
this direct method and would go even a little farther than you on
p. 4 of your manuscript. I am not of the opinion that the later
method is the more progressive one.
\cite{SchuranWeyl:Hs91:739}\footnote{Schur's (undated) letter is
an answer to a letter by Weyl, which is not preserved. The
discussion relates well to \cite{Weyl:79}. The only point I
cannot identify is the  the reference to the  remark ``\ldots  on
p. 4 of your manuscript \ldots''.}

\end{quote}

Any  representation of $S_f$ is characterized by a
character $\chi $, i.e., the  complex valued function on  $S_f$, defined by the trace of the corresponding represention matrices. For an irreducible representation it is known that $(\chi , \chi ) = 1$, with respect to the  scalar product in the  function space on $S_f$.  In the sequel we shall use the notation $\rho _V (\chi )$ for
the irreducible representation of $GL(V)$ in $\bigotimes^f V $,
corresponding to $\chi$ by this correspondence and Weyl's theorem.\footnote{Weyl's notation
for our $\rho_V (\chi )$ was $\Lambda _n (\chi )$, where $n = dim
\, V$. } Weyl  considered a spin-extension of the underlying
vector space of 1-particle states, $V$ ($ dim \, V = n$), in the
sense  of Pauli wave functions, 
\beq  \label{otimes} W:= V   \otimes \C^2
 \; , \, \;\; dim \, W = 2n \,  . \eeq 
In the case  of
an $f$-electron system one has to study the irreducible
components of the  operation of  $GL(V)$ induced on the
antisymmetric part of the tensor product, $\bigwedge ^n W$. The
decomposition of  $\bigwedge ^n W$ according to Weyl's main theorem
leads to multiplicities $m_{\chi }$ for the irreducible
representations of type $ \rho _{W} (  \chi )$, such that \beq
\label{decomposition}  { \bigwedge} ^n W = \bigoplus m_{\chi } \,
\rho _{W} (  \chi  ) \,  \eeq For the   calculation of  the
multiplicities $m_{\chi }$ Weyl  established  a kind of ``duality'' (Weyl's terminology) among  the  representations  of the symmetric group.

To any representation $\rho_{U}$ of  $S_f$ in a vector space $U$ there is an induced representation $\rho_U^{\ast}$ on the dual  space $U^{\ast}$. By contextual reasons, Weyl  modified the sign of this induced operation on $U^{\ast}$ by the signum function.\footnote{If $\rho_U$  corresponds to a character $\chi$,   Weyl defined the  {\em dual representation}  $\chi^{\ast}$  as the representation of $S_f$  given by
\[  \sigma \mapsto signum (\sigma) \rho_U^{\ast}(\sigma) \quad \mbox{ \cite[187]{Weyl:79}}. \] }
Then he could use the apparatus of character formulae and found a striking {\em reciprocity  relation} (Weyl's terminology) between the  multiplicity of an irreducible representation of the symmetric group and the dimension of its dual representation:

\begin{theorem}
 The multiplicities $m_{\chi }$   in (\ref{decomposition}) are equal
to the dimensions of the corresponding  dual representations
$\chi ^{\ast}$, 
 \[  m_{\chi } = dim \,  \chi  ^{\ast}   \; ,  \]
\cite[187]{Weyl:79},\cite[352]{Weyl:GQMEnglish}.
\end{theorem}

A  direct consequence was that $m_{\chi } =0$, if the
Young diagram corresponding to $\chi $ has more than 2  columns.\footnote{The signum factor in Weyl's definition  of the dual representation implies $ dim \,   \chi  ^{\ast} = 0$ for dual representations with more than 2 rows.  The  Young diagram of the representation  in the dual space $U^{\ast}$  is obtained from the diagram in $U$ by transposition. Thus only representations with Young  diagrams  of 1 or 2 columns have non-vanishing multiplicities in the decomposition of the alternating product
(\ref{decomposition}) \cite[350, 352, 370]{Weyl:GQMEnglish}.}
From a pragmatic point of view, this result stated the same  condition for the existence of an antisymmetric spin extension as the one given by Wigner and von Neumann in terms of the partition $(\lambda)$ (equation (\ref{2-2-1-1})).
But Weyl considered this insight as more than just a calculational
tool. For  him it established a kind of {\em reciprocity law} of
undoubtedly {\em material importance}.
\begin{quote}
The modification, which is brought about  by the existence of spin under neglection of its dynamical effects and by the Pauli exclusion principle, consists  in nothing more than in  {\em  a transformation of the  multiplicity  of the term system corresponding to $\chi$ from }  [$m_{\chi }$] {\em   into} [$dim \,  \chi  ^{\ast}  $].  \ldots 
The dynamical effect of spin resolves these {\em multipletts} in as many components, as given by its multiplicity  [$dim \,  \chi  ^{\ast}  $]; moreover it induces  weak intercombinations between the different classes of terms. [Notation adapted to ours, emphasis in original, E.S.] 
\cite[188]{Weyl:79}
\end{quote}

\subsubsection*{Spin coupling in  general exchange molecules }
Weyl  even  extended the reciprocity theorem  to
 a more general case,
$W    = V'  \otimes V''$.  
At first glance, this generalization may look like a pure mathematicians game,  without connections to the physical context,  but Weyl was highly interested in its application to molecular bonds. 

He considered  two atoms $A$ and $B$  with  electron numbers $\nu '$ and $\nu ''$ and symmetry types given  by the   irreducible representations ${\cal
G}_{\chi ' }$ ,  ${\cal G}_{\chi '' }$ (with  characters $\chi '$ and $\chi ''$ ---  Weyl's notation). If they form  a molecule, the bond would be described
by (collective) states of the combined electron system in the
tensor product. The mathematically elementary states would then
be characterized by the irreducible representations in the
product.  Weyl  generalized Heitler's and London's theory from exchange molecules  with electron pairs  to the many
 ($\nu = \nu '  + \nu '' $) electron case. His generalized reciprocity theorem
(Weyl's terminology) contained the clue for analyzing  the possible
bonding constellation of higher atoms.

 In one of his presentations of the result to  a wider audience, a
published version of  talks he gave during his journey through
the United States in late 1928 and early 1929, he  explained his
basic idea:
\begin{quote}
This reciprocity law governs the fundamental chemical problem of
combining two atoms to obtain a molecule \ldots . The molecule
which is obtained by combining the two atoms will be in one of
the symmetry states $\zeta $ whose corresponding ${\cal G}_{\zeta
}$ [Weyl's symbol for an irreducible representation of the full
permutation group of all $\nu = \nu ' + \nu ''$ electrons with character $\zeta$, E.S.]
appears in  ${\cal G}_{\chi '} \times {\cal G}_{\chi '' } $ and
the calculation of the associated energy is accomplished with the
aid of these characteristics [characters, E.S.]. These
circumstances which cannot be represented by a spacial (sic!)
picture, constitute the basis for the understanding of the
homopolar bond, the attraction (or repulsion) existing between
neutral atoms \ldots \cite[290f.]{Weyl:79}
\end{quote}

With respect to the strong conceptual relationship between mathematics and physics, 
these words may appear similar to  those  Weyl had written a
decade earlier, in the years between 1918 and 1920 when he
pursued his program of  a geometrically unified field theory. But 
 during the 1920s  Weyl had become  much more sensitive to empirical questions. At the end of the decade he had  the impression that ground was 
touched  in the formerly fathomless search for a  mathematization of the
basic structures of matter. This new viewpoint seemed
incompatible with the earlier hopes  for a unified field theory
of matter in terms of classical fields, which Weyl now considered to be 
 illusionary.\footnote{Compare \cite{Scholz:WeylMatter}.}   The  role   played in his earlier work  in general relativity and  unified field by generalized differential
geometric structures  was now taken over
by   group representations in Hilbert spaces (``unitary
geometry'') and the quantum theory of atoms and their bonds.

While in the early 1920s he still thought in terms of   {\em a-priori} structures supported by strong methodological and ontological speculations,
he now only spoke of an  ``appropriate language'' for the
expression of the natural  ``laws''.
\begin{quote}
The connections between mathematical theory and physical
application which are revealed in the work of Wigner, v. Neumann,
Heitler, London and the speaker is here closer and more complete
than in almost any other field. The theory of groups is the
appropriate language for the expression of the general
qualitative laws which obtain in the atomic world. (ibid.)
\end{quote}

In winter 1928/29  Weyl  used a journey   to the US to bring the
gospel of group theory to the  scientifically rising country. He gave
lectures at Princeton and Berkeley on ``Application of group
theory to quantum mechanics''  \cite{Weyl:Hs91a:51}, and
published three articles on the topic in North-American journals
\cite{Weyl:79,Weyl:86,Weyl:87}.\footnote{\cite{Weyl:79} was
published in German in  the {\em Annals of Mathematics}.} After
his  move from Z\"urich to G\"ottingen in early 1930, he took
part in the seminar on the structure of matter, which went back to the 
Hilbert tradition and was now run by  Born. He was thus led to a
further elaboration of his method \cite{Weyl:90,Weyl:91}. The
second of these notes contained an analysis of determinantal
methods used by W. Heitler and G. Rumer in their common work
presented in the seminar.\footnote{ \cite{Heitler/Rumer:1931}}

Building on his previous analysis, Weyl showed how to express the
spin states of an $m$-electron system formed from the shells of $k$ atoms, with
$m_1, \ldots, m_k$ valence electrons each  ($m = \sum_1^k
m_j$), and the condition that
 $m_0$ valences remained free. Admissible spin coupling constellations of the valence
electrons could be constructed  from alternating products  of the
eigenfunctions of pairs of electrons from different atoms. After
assigning variables $x_1, \ldots, x_k$ to each atom and $x_0$ to represent empty valences,  Weyl developed a  method to calculate molecular bond energies. The method relied on the
{\em first fundamental theorem of invariant theory} according to
which it is possible to express  the invariants of any set of  vectors
 $ \{ x_0, \ldots,  x_k \} \subset \C^2$ under the operation of 
  $SL_2(\C )$
by integer polynomials in  the ``fundamental invariants''  $z_{i,j}$ derived from the vectors  by determinants
\[ z_{i,j} := det (x_i, x_j) \qquad 0 \leq i, j \leq k \; . \]
Weyl used the abbreviated notation $ z= [x,y] $ (the {\em fundamental binary invariant}), for any two vectors $x$ and $y$.

 According to Weyl, a  ``pure valence state'' was  characterized 
by a monomial  of total order $m$ and order $m_j$ in each
component $x_j$  ($0 \leq j \leq k$), formed from  binary
invariants $[x,y]$.\footnote{The totality of  pure valence states is
not algebraically independent, but obeys  a relation, given by
the ``second fundamental theorem of invariant theory''. }
Eigenstates of the  molecule would not be pure valence states but  {\em superpositions} of them, which are eigenstates of the Hamilton $H_p$ operator of the  bound and spin perturbed
system, 
\[ H_p = H_0 + \sum H_{\alpha  \beta }  \; ,\]  linearized in  terms  due to the exchange (transposition) of any two of the valence electrons. Here $H_0$ denotes the Hamilton operator of the electron system without spin coupling.
 Weyl developed a method for a  calculation of the perturbation term $H_p - H_0$, if the  exchange energies  $W_{\alpha \beta }$ between  two valence electrons ($1\leq \alpha  \leq m_i, \; 1 \leq \beta \leq m_j$) of two atoms with index $i$ and $j$  could be calculated \cite[323f.]{Weyl:91}.
 The critical
point for applications of the method was then the 
calculation of all the ``exchange energies'' involved. It presupposed the solution of a generalizated version of Heitler's and London's problem for  electron pairs. Moreover, the whole method could be physically relevant only  for molecules for which the exchange energy contributes essentially to the total bond energy. Molecules with large $H_0$, with respect to the spin perturbation, could be analysed  just as well by studying only the 
Schr\"odinger wave component of their Pauli spinors.\footnote{These are molecules in which the geometry of ``molecular orbits'' of valence electrons and the Coulomb potential are the  essential  determinants of  the bond energy.}

From a theoretical perspective,the structure of
the procedure was very satisfying. Weyl argued that, by assigning
formally a ``valence dash'' (between atom $x$ and $y$) to each binary invariant of type $[x,y]$, one arrived at 
 graphs for pure valence states, which were in striking agreement with an old proposal by
J.J. Sylvester.  In 1878, Sylvester had proposed,  in a purely
speculative approach, to  express chemical valence relations  by
binary invariants.  Formally his proposal coincided 
with the algebraic core of Weyl's  construction. Now Sylvester's procedure could be understood as an expression of an
algebraic structure underlying  the determination of bound 
states in the new quantum mechanical theory of valence bonds. No wonder that Weyl and 
Heitler  were fond of the new quantum chemical underpinning of
Sylvester's  speculative method.\footnote{For a more detailed
discussion see \cite{Parshall:Chemistry} and \cite[section 3.1,
163--177]{Karachalios:Diss}.} 

There remained, of course, several problems. 
The practical usefulness
of the method could be tested only if the 
exchange energies of single electron pairs could somehow  be calculated. Even then it remained to be seen, whether the result would be in agreement with empirical chemical knowledge. In
his first publication,  Weyl  only indicated  the 
general method
\cite[323f.]{Weyl:91}.\footnote{A graphical method for the
construction of a basis of invariants, based on an idea of G.
Rumer,  was  written down by Rumer, Teller and Weyl in
\cite{Weyl:97}.}
 In the 1930s he continued with the calculation of  examples. That is 
 shown by notes in his 
  {\em Nachlass} \cite{Weyl:Hs91a:32} and by  remarks in a
new appendix written for   \cite{Weyl:PMNEnglish}. 

But the method was never adopted in the chemical community.  Most of the molecules of organic chemistry turned out to be different from the bonding class which Heitler had called exchange molecules, even  in Weyl's generalization. 
 During the years, chemists found overwhelming evidence that their  models of molecular orbits, in which the spatial distribution of the Schr\"odinger part of the wave function contributed decisively to the binding energy and sufficed in most cases  to solve their problems.  Moreover, the  method  of molecular orbits was closer to the imagination of the chemists and  its mathematics was  easier to handle for them.
The more structural  method of exchange energies of spin coupling
remained marginal for the practice of physical chemistry, even in the extended and refined form which Weyl had started to develop and to present as a methodological tool to the community of physicists and physical chemists .

\section*{7. Outlook}

In spite of its surprising theoretical achievements,
 the rise  of  groups in quantum mechanics was far from a straight forward
 story. With  its first  successes at the turn to the 1930s,
there arose sceptical reservation, criticism,  and even strong
counterforces to the spread of group theoretic methods in the new
field of theoretical physics. Such criticism was not always meant
as a real opposition to the modernizing tendency; sometimes it
was just an    expression of uneasiness with  the new algebraic methods.  
Soon after  Pauli moved from Hamburg
to Z\"urich as the successor of Debye,   in April 1928, Ehrenfest
asked him for help in the difficult the new matter. Pauli
was well-known for his ability to absorb new mathematics with
ease and to adapt it to the necessities of theoretical physics.
Moreover, in his last year at Hamburg he had participated in a
lecture course on algebra and group theory given by Emil Artin. After his arrival in Z\"urich in early 1928,  
he stood once again in close communication with Weyl like in the early 1920s.\footnote{\cite{Meyenn:Pauli_Symmetries},
\cite[472]{Mehra/Rechenberg:VI}. A couple of weeks after his
arrival Z\"urich, Pauli wrote in a letter to N. Bohr : ``I have
now learned so much erudite group theory from Weyl that I am
really able to understand the papers of Wigner and Heitler''
\cite{PaulianBohr:201a}. Moreover, he read and commented page
proofs of Weyl's  GQM in early summer 1928
\cite[402]{PaulianWeyl:2183}.}

\subsubsection*{Group pest}
In September 1928, Ehrenfest turned to Pauli and asked  for help in
understanding the ``terribly many papers on the group-pest
(Gruppenpest)'', of which he ``could not read any one beyond the
first page'', as he wrote to Pauli on September 22,
1928.\footnote{Quoted from \cite[473]{Mehra/Rechenberg:VI}.} In  parts of the --- still small --- 
community, this
word   became the catchword for opposition to the use of 
group theoretic methods  in quantum mechanics. Apparently  Ehrenfest 
 unwillingly contributed a verbal battle sign to  the emerging anti-group
camp. For him the word  expressed nothing more
than uneasiness about the rising challenges of the new
mathematical methods in  theoretical physics. He was not at all  opposed  by principle to the new
tendencies. On the contrary, he supported its development
actively. On his initiative,  B. L. van der Waerden started to
develop his calculus of spinor representations of the Lorentz
group \cite{Waerden:Spinoren}; and one of his later doctoral
students, H. Casimir, started to do research work on quantum
mechanics, very much  influenced by Weyl's book. As has been
discussed on other occasions,\footnote{\cite{Meyenn:Casimir},
\cite[512--514]{Mehra/Rechenberg:VI}, \cite{Hawkins:LieGroups}.}
  Casimir  finally even contributed
to the refinement of  representation theory itself, by
proposing  an idea for a purely algebraic proof of the full
reducibility of representations of Lie groups, derived from his
research on the problem of rotation in quantum
mechanics.

Real and strong opposition to the group theoretic approach to
quantum mechanics came from another camp led by John Slater, who
showed that already  traditional  algebraic tools   were highly
effective in the calculation of the energy of higher atoms and
binding energies of molecules \cite{Slater:1929}. Slater's
background in a more pragmatic tradition of  theoretical physics
in the United States surely played a role for his strong
rejection of the more theoretically minded approaches like
representation theory  \cite{Schweber:Slater}.\footnote{See also
\cite{Skuli:Diss}, \cite[499ff.]{Mehra/Rechenberg:VI} and  for a
broader comparative discussion of German and American physical
chemists  of the first generation  \cite{Gavroglu/Simoes}.}

Slater's success in developing determinant methods for quantum
mechanical calculations  found immediate acceptance among leading
protagonists of the G\"ottingen milieu. Shortly before  Weyl
decided to come back to G\"ottingen as the successor to David
Hilbert, Max Born warned him, in an otherwise very friendly
welcome letter, that he supported the  ``attempt to throw group
theory out of the theory of atomic and molecular structures, as
far as possible'' \cite{BornanWeyl:Hs91:490}. At that time, Born
was  close to finishing an article in which he attempted to get
rid of group theoretic methods in the theory of chemical bonds
\cite{Born:1930}. He even was proud of having convinced Heitler,
after the latter's arrival at G\"ottingen as Born's assistant, to
give up the idea that  group theoretic considerations might play
an important role in studies of molecular
bonds.\footnote{\cite{BornanWeyl:Hs91:490}} This perspective
resulted in a common article by W. Heitler and G. Rumer on
chemical bonds, which only used ``traditional'' algebraic methods
along the line of Slater and Born
\cite{Heitler/Rumer:1931}.\footnote{The article  was written
after  Weyl had arrived at G\"ottingen, and after a discussion of
the method in the common seminar on the structure of matter.}
On the other hand, group theoretic methods in physics and quantum chemistry 
continued to be a topic for  lecture courses at the G\"ottingen mathematical institute.\footnote{W. Heitler gave a course on this subject in winter semester 1929/30 at the mathematical institute \cite{Heitler:VorlesungWS29/30}. He  concentrated on the subject matter 
of  Wigner's and von Neumann's theory. Only in the last chapter he gave a short introduction to the theory of molecular bonds. I owe Martina Schneider the information on this course.}

\subsubsection*{Weyl at G\"ottingen}
In the meantime, in May 1930, Weyl had  accepted the call to
G\"ottingen and started to teach there in
winter semester of the same year. That gave him a splendid
occasion for critical exchanges  and   collaboration with  Born,
 Heitler,   Rumer, and  Teller on   group theoretical
methods in the nascent quantum chemical context. Although  Born  had been highly sceptical of the method earlier on,  he  gave critical support to  the
enterprise after Weyl moved to G\"ottingen, in his own way. This
exchange of ideas with the  theoretical physicists around Born in
the common G\"ottingen seminar led Weyl to a more detailed
elaboration of his use of symmetry operators in the $n$-fold
tensor space of electron states for the characterization of
molecular bond states and the establishment of the link to binary
invariants \cite{Weyl:90,Weyl:91}. In a
subsequent review article  on the quantum theory of molecular
bonds in the {\em Ergebnisse der exakten Naturwissenschaften},
Born finally rephrased those results  of Weyl's investigation
which seemed of importance to him for physicists and physical
chemists. In the introduction to his article he frankly declared
that the proofs of Weyl's results could not be rephrased under
``complete avoidance of the `group pest' which Slater and
the author [Born] had intended''. He therefore restricted the
presentation to formulas and rules, without proofs, such that the
results could be understood by physicists and chemists without
being forced to read ``the difficult works of Frobenius and Schur
on the representation theory of groups'', as he wrote in his
introduction \cite[390]{Born:1931}.

All in all, the first wave of rapid development of group
theoretical methods in quantum mechanics ran into the opposition of
a strong,  multi-faceted,  anti-group camp; or, at least, it
had to face pragmatic  scepticism among physicists and theoretical
chemists at the turn to the 1930s.

On the other hand, new forces joined the party of mathematical
contributors to representation theoretic methods for mathematical
physics. Most important,  from the side of young mathematicians,
was Bartel Leendert van der Waerden  who entered this scene   with his
spinor  paper written with the explicit goal of serving the
physics community \cite{Waerden:Spinoren}.\footnote{More  details
will be discussed in \cite{Martina:Diss}.}
 In personal communications with Weyl he also contributed critical
 remarks to the  understanding of   algebraic structures underlying spin
coupling.  Van der Waerden  criticised Weyl's approach  from the
viewpoint of a  young ``modern'', i.e., structurally oriented,
algebraist. In  a letter from April 4, 1930, he argued that  in
Weyl's derivation of the ``reciprocity theorem''  it was
unnecessary  to build upon  the ``inessential property that $\pi$
[Weyl's symbol for the  permutation group, E.S.] is a permutation
group''. Obviously  he  abhorred the ``multitude of indices'' used
by Weyl and claimed that one could do without them in this
investigation\cite{WaerdenanWeyl:Hs91:784}. After some exchanges
of letters, of which only the van der Waerden part is preserved,
he argued that the result was essentially a question in the
representation theory of algebras. According to van der Waerden's
analysis, Weyl's result depended  essentially on the fact that a 
 matrix algebra ${\cal A}$ induced from the operation of the group algebra
$\C [S_f]$ on $\bigotimes^f V$  commutes with a completely reducible representation of the
general linear group $GL(V)$  on the tensor product
\cite{WaerdenanWeyl:Hs91:785}.\footnote{See also
\cite{Waerden:1930}.}
 It seems that Weyl was not completely convinced that such a level of structural abstraction suited his purpose. He rather    insisted on the  use  of
 the ``multitude of indices'',  because they were essential for
the context of modelling the combined electron systems of two
atoms in a molecule. Nevertheless he accepted the proposal to
straighten the derivation of the reciprocity theorem
\cite[310]{Weyl:91}.

In this sense, the interaction between physicists and mathematicians  close to the
 G\"ottingen and Z\"urich milieu seemed to be a
 a splendid scientific environment for a further consolidation of
 group theoretic methods in physics and chemistry at the turn to
the 1930s. In the next couple of years, the triad of now
classical text books on the use of group theory in quantum
mechanics appeared \cite{Wigner:Gruppen}, the second edition of Weyl's GQM
and  its English translation by H.P.
Robertson \cite{Weyl:GQMEnglish}, and \cite{Waerden:GQM}. These
books broadened the  basis for an extension of the approach,
invited scepticists to take an own look at the question, and
enabled newcomers from different backgrounds to join the
enterprise.

\subsubsection*{From an intermediate period \ldots}
As we know, and most of the participants sensed  well, the social
stability of this milieu stood on shaky ground. Only little
later, with the  Nazi's rise to power,  the G\"ottingen
mathematical science group was dismantled. As one of the
consequences,    the closely knit interaction between pragmatic
sceptics with respect to the group theoretic method, close to
Born, and the group of active  protagonists like Weyl, van der
Waerden, Heisenberg, Wigner and von Neumann, which was easily organized around
G\"ottingen, was interrupted. Although several of the
protagonists of the first wave continued to elaborate and to teach or propagate the new
method, no great gains in terms of
broader acceptance could be made during the next two decades .

Weyl continued to argue for the use of  the new method, in
particular in the context of chemical bonds, in publications,
talks and lecture courses. But he was very well aware of the
reservations of  the practitioners of the field felt in rageard to 
his proposals of using invariant theory for the  characterization of  
bond states, and he accepted it. In an undated manuscript
of a talk given in the second part of the 1930s, Weyl remarked
that the development in the field had not been ``very favorable
to the scheme'' which he had laid out. The recent report
\cite{vanVleck} had nearly passed it over ``in silence''. He
realistically added that in his exposition he even  intended to ``clearly
indicate the boundaries of applicability for our scheme''
\cite[2]{Weyl:Hs91a:32}.

Finally he  concentrated his research and publication efforts
on the mathematical foundation of the theory. In joint work with
Richard Brauer he developed a global characterization of spin
representations in any dimension (and of arbitrary signature) by
Clifford algebras \cite{Weyl/Brauer}.\footnote{E. Cartan had
discussed spinor representations on the infinitesimal level
already in 1913;  here the integral (global) perspective stood in
the center. } All this
 culminated in his book on {\em The Classical Groups} \cite{Weyl:ClassGroups}.
That was no disillusioned
withdrawal to pure mathematics. It rather was an   expression
of a realistic evaluation of the actual situation in the field of
application. Even though Weyl's calculation of binary invariants
did not enter the core of the theory of chemical bonds, his
invariant theoretical analysis of spin constellations turned out, in the
long run, to be an important contribution to the study of   spin-coupling, which has recently started to attract 
new interest from the point of view of ``entangled'' systems. The introduction
of binary invariants into the study of coupled systems of
electrons in the late 1920s and the following decade, may turn
out to be another prelude to the development of a symbolic game
with long lasting importance in a shifted context of
application.\footnote{This ``game'' has recently gained new interest from the point
of view of {\em quantum computing}. In  this new context  the question of  energy contributions, which
hindered Weyl's proposals from becoming  important in 
 quantum chemistry,  are subordinate. I owe  the hint to the
connection of Weyl's work with these recent developments to P. Littelmann.}

During the decades of slow maturation, it  was mainly due to Werner Heisenberg's anticipatory guess of {\em isospin} $SU_2$ as a symmetry underlying the nuclear interactions \cite{Heisenberg:Isospin} and to  Eugene Wigner's continuing work and
insistence on the importance of the group theoretic approach for
fundamental physics, that this research tradition in mathematical physics was never completely
interrupted.\footnote{Cf. \cite{Rasche:Isospin} and \cite[265f.]{Mackey:Wigner}.}   Most important  for relativistic  quantum physics was Wigner's fundamental
work on the representation theory of the Poincar\'e group
\cite{Wigner:39}.

\subsubsection*{\ldots to a second wave of groups in  quantum physics}
With the exception of  such ``heroic'' but for a long time
relatively isolated contributions, it needed a new generation of physicists
and a diversification of problems and another problem shift in
quantum physics,   before group theory was stepwise integrated
into the core of quantum physics. Faced with the rise in
complexity of problems of nuclear spectroscopy, G. Racah brought
group theoretic methods closer to the ordinary problem solving
practice of spectroscopists
\cite{Racah:1942-49}.\footnote{\cite[269]{Mackey:Wigner}} Finally  the proliferation of new ``elementary particles''
between 1950 and the 1970s gave material and motivation  to look  for group theoretical classifications of object
structures and the corresponding internal symmetries of  interactions. Thus we can see a second wave in the use  of group theoretical methods in quantum physics
during the 1950s to the 1970/80s. In this changed context, the
two books of  the  above mentioned triad, which formerly were only
available in German, were translated into English,
\cite{Wigner:Groups} and \cite{Waerden:GQMEnglish}. Mathematicians
of the next   generation, among them G. Mackey and I.E. Segal,
continued to contribute, from the side of 
mathematics, to the research tradition begun at the end of the
1920s.

In this second wave of research, simple anticipatory ideas had to be differentiated and  different strands of using groups
in quantum physics grew together:
\begin{itemize}
\item[---]  weight systems of representations were turned into a tool for understanding ``multipletts'' of  basic states of matter,  generalizing  the multipletts of spectral terms of the 1920s,
\item[---] isospin was first enriched (``eightfold way'', $SU_3$) and then transformed into two different forms (weak isospin, $SU_2$, and the ``chromo-symmetry'' of strong interactions, $SU_3$),  the basic symmetries of  particle physics of the late 20th century,
\item[---] conservation laws became generally considered as founded upon underlying dynamical symmetries, 
\item[---]  the study of  infinitesimal
symmetries became standardized in the form of  (generalized and non-abelian) gauge fields or, equivalently,  connections in fibre bundles.
\end{itemize}

 Groups, their representations, corresponding conserved quantities, and the use
 of  gauge structures were finally broadly accepted. They were  used
as  an important ingredient of  the  mathematical  forms functioning
as  a  symbolic relative a priori in which theoretical
physicists of  the late 20th century were able to mold an
impressive part of the   experimental
knowledge of fundamental physics. 
At the end of the second wave, group theoretical  methods were  well   integrated into the mainstream of mathematical
physics.  Although  at the end of the century  the gap  between general relativity and quantum physics continued  to be  wide open,  groups and their representations have turned into useful tools and provide conceptually convincing forms for the 
construction of symbolic models of material processes  in both domains.

\small
\section*{Acknowledgements: }
This paper is the result of an extended period of study and
discussions, interrupted by other activities. I am indebted to the support of the {\em
Handschriftenabteilung der ETH Bibliothek Z\"urich} and the {\em
Staatsbibliothek Berlin} for access to the archival
material. Dr. Michael Weyl  kindly permitted to use  his father's  {\em Nachlass} at
Z\"urich. The {\em Max Planck Institute for History of Science}, Berlin,  gave me the
opportunity for first  discussions on the history of
quantum mechanics with colleagues from the history of physics in
November and December 1998. I got important hints on diverse
aspects of the topic in conversations with  many colleagues,
most importanly, although now years ago,  with Peter Slodowy. More recently G\"unther Rasche, Jim  Ritter, David Rowe, Tilman
Sauer, Martina Schneider, Sk\'uli Sigurdsson, Urs Stammbach, Christian Wenzel and an anonymous referee of the manuscript contributed by their comments. Finally the  interest and
patient support by Jeremy Gray, including a final revision of language, was decisive for writing down this article.


\end{document}